\documentclass[11pt]{article}
\usepackage{fullpage,amsthm,enumitem}
\usepackage{amsmath,amsfonts,amssymb,graphicx} 

\usepackage{mathpazo} 
\usepackage[dvipsnames]{xcolor} 
\usepackage[hyperfootnotes=false]{hyperref}
\usepackage[title]{appendix}
\usepackage{tensor}

\usepackage[backend=biber,
  style=alphabetic,
  giveninits=true,
  maxnames=10,
  maxalphanames=4,
  doi=false,
  url=false,
  isbn=false,
  eprint=false,
]{biblatex}
\renewbibmacro{in:}{}
\DeclareFieldFormat
[article,inbook,incollection,inproceedings,patent,thesis,unpublished]
  {title}{#1\isdot}

\newtheorem{thm}{Theorem}
\newtheorem{cor}[thm]{Corollary}
\newtheorem{lem}[thm]{Lemma}
\newtheorem{prop}[thm]{Proposition}
\newtheorem{defn}[thm]{Definition}
\newtheorem{ex}[thm]{Example}
\newtheorem{rem}[thm]{Remark}

\newcommand{\N}{\mathbb{N}}
\newcommand{\Z}{\mathbb{Z}}
\newcommand{\R}{\mathbb{R}}
\newcommand{\C}{\mathbb{C}}

\renewcommand{\H}{\mathcal{H}}

\newcommand{\A}{\mathcal{A}}

\newcommand{\mB}{\mathcal{B}}

\DeclareMathOperator{\Tr}{Tr}
\DeclareMathOperator{\dom}{dom}

\newcommand{\op}{\operatorname{op}}
\newcommand{\OP}{\operatorname{OP}}
\def\multinom#1#2{\ensuremath{\left(\kern-.3em\left(\genfrac{}{}{0pt}{}{#1}{#2}\right)\kern-.3em\right)}}

\title{
\Large 
Multiple operator integrals, pseudodifferential calculus, and asymptotic expansions
}
\date{\today}

\author{Eva-Maria Hekkelman, Edward McDonald, Teun D. H. van Nuland}

\begin{document}

\maketitle{}
\begin{abstract}
    We push the definition of multiple operator integrals (MOIs) into the realm of unbounded operators, using the pseudodifferential calculus from the works of Connes and Moscovici, Higson, and Guillemin. This in particular provides a natural language for operator integrals in noncommutative geometry. For this purpose, we develop a functional calculus for these pseudodifferential operators. 
    To illustrate the power of this framework, we provide a pertubative expansion of the spectral action for regular $s$-summable spectral triples $(\mathcal{A}, \mathcal{H}, D)$, and an asymptotic expansion of $\mathrm{Tr}(P e^{-t(D+V)^2})$ as $t \downarrow 0$, where $P$ and $V$ belong to the algebra generated by $\mathcal{A}$ and $D$, and $V$ is bounded and self-adjoint.
\end{abstract}

\section{Introduction}
\subsection{Operator Integrals}
\footnote{2010 \emph{Mathematics Subject Classification}:47A55; 47A60, 47G30, 58B34, 47F10.}
\footnote{\emph{Key words and phrases}: Multiple Operator Integrals, Noncommutative Geometry}
Operator integrals appear in various areas of noncommutative geometry (NCG)~\cite{Connes1994}, one of the most prominent examples of which is the JLO cocycle~\cite{JLO1988}
\[
\int_{\Sigma^n} \Tr(\eta a^0 e^{-t_0 D^2} [D,a_1] e^{-t_1D^2} \cdots  [D,a_n]e^{-t_nD^2})dt .
\]
The expansions and other identities derived for these operator integrals are often very similar in nature, and it would therefore be useful to have a framework to systematize these techniques. The language tailored to provide such a framework is that of multiple operator integrals (MOIs)~\cite{ACDS, Peller2006, Peller2016, SkripkaTomskova2019}. However, the existing literature on MOIs is not equipped to make sense of multiple operator integrals as pseudodifferential operators in the style of Connes--Moscovici, Higson and Guillemin~\cite{ConnesMoscovici1995, Guillemin1985, Higson2003}, which would be necessary for advanced applications in NCG.

We therefore first generalise the MOI construction of \cite{ACDS, Peller2006}, realising the operator integral as an abstract pseudodifferential operator on the condition that the operator integrand is itself an abstract pseudodifferential operator. Then we develop a functional calculus for abstractly elliptic pseudodifferential operators, which completes the technical framework. To ensure that divided differences are valid symbols for our MOIs, we prove a novel result on integral decompositions of divided differences $f^{[n]}$  using almost analytic extensions. Finally, we derive various identities for these MOIs that are often used in the NCG literature, and conclude by showing that these identities can be used to provide a perturbative expansion of the spectral action and prove a new result on the existence of asymptotic expansions in regular $s$-summable spectral triples $(\A, \H, D)$ of
\[
\Tr(P e^{-t(D+V)^2})
\]
as $t\downarrow 0$, where $V$ is some self-adjoint bounded element in the algebra generated by $D$ and $\A$, and $P$ is an arbitrary element in this algebra. The assumptions needed to deduce the existence of these asymptotic expansions was an open problem, communicated to the authors by Bruno Iochum.

Multiple operator integrals find their origin in the work of Daletski\u{\i}--Kre\u{\i}n~\cite{DaletskiiKrein1956}, Birman--Solomyak~\cite{BirmanSolomyak1,BirmanSolomyak2,BirmanSolomyak3}, Pavlov~\cite{Pavlov1971]} and Sten'kin~\cite{SolomyakStenkin1971,Stenkin1977}. In the theory that we adapt from~\cite{ACDS,Peller2006}, classical MOIs have a symbol $\phi: \R^{n+1} \to \C$ which can be written as 
\begin{equation}\label{eq:BSrep}
\phi(\lambda_0, \ldots, \lambda_n) = \int_\Omega a_0(\lambda_0, \omega) \cdots a_n(\lambda_n, \omega) d\nu(\omega),
\end{equation}
for some finite measure space $\Omega$ and bounded measurable functions $a_i: \R \times \Omega \to \C$. A MOI with such a symbol is the well-defined operator 
\begin{equation}
T_\phi^{H_0, \ldots, H_n}(X_1, \ldots, X_n) := \int_\Omega a_0(H_0, \omega) X_1 a_1(H_1, \omega) \cdots X_n a_n(H_n ,\omega) d\nu(\omega),
\end{equation}
where $H_i$ are closed densely defined self-adjoint operators on a separable Hilbert space, and $X_i$ are bounded operators. For $n=0$ this is simply functional calculus, but for $n > 0$ MOIs are multilinear maps
\begin{align*}
    T_\phi^{H_0, \ldots, H_n}&: \underbrace{B(\H)\times \cdots \times B(\H)}_{n} \to B(\H).
\end{align*}
Here, and throughout the paper, $B(\H)$ denotes the algebra of all bounded linear operators on a separable Hilbert space $\H$. The notation $T_\phi^{H_0, \ldots, H_n}$ is justified since this multilinear map does not depend on the chosen representation~\eqref{eq:BSrep} (see \cite[Lemma~4.3]{ACDS}). To ease notation, if $H_0 = \cdots = H_n =: H$ we simply write $T_\phi^{H} := T_\phi^{H_0, \ldots, H_n}$.

The quintessential symbol for MOIs is a \textit{divided difference}, namely, for $f:\R \to \C$ regular enough we can recursively define
\[
f^{[0]}(\lambda_0):=f(\lambda_0), \quad f^{[n]}(\lambda_0, \ldots, \lambda_n): = \frac{f^{[n-1]}(\lambda_1, \ldots, \lambda_{n}) - f^{[n-1]}(\lambda_0, \ldots, \lambda_{n-1})}{\lambda_0 - \lambda_{n}}, \quad \lambda_0 \not = \lambda_{n},
\]
replacing the above by a limit in case $\lambda_0 = \lambda_n$, so that
\[
f^{[n]}(\lambda, \ldots, \lambda) = \frac{1}{n!}f^{(n)}(\lambda).
\]
For a nice enough function $f$ the divided difference $f^{[n]}:\R^{n+1}\to \C$ admits a decomposition of the required form, and it then holds that
\[
\frac{1}{n!}\frac{d^n}{dt^n}\bigg|_{t=0} f(H+tV) = T_{f^{[n]}}^{H}(V, \ldots, V) \in B(\H).
\]
For commuting operators, this formula reduces to the chain rule.

MOIs with divided differences as their symbols appear in NCG in the context of spectral flow~\cite{AzamovCarey2007}, spectral shift~\cite{ACDS}, the spectral action~\cite{Skripka2014,Skripka2018, Suijlekom2011}, the heat trace expansion~\cite{lesch2017divided,NulandSukochev2023}, and cyclic cocycles~\cite{CPRS1, CPRS2, Liu2022, vNvS21a}. The JLO cocycle mentioned earlier can explicitly be written as 
\[
\int_{\Sigma^n} \Tr(\eta a^0 e^{-t_0 D^2} [D,a_1] e^{-t_1D^2} \cdots  [D,a_n]e^{-t_nD^2})dt = \Tr(\eta a_0 T_{f^{[n]}}^{D^2}([D,a_1], \ldots, [D,a_n])),
\]
where $f(x) = \exp(-x)$. The properties and expansions of these cocycles are derived from a small list of MOI-identities, which have also proved instrumental in \cite{ACDS,Liu2022,Peller2006,PSS2013,PS,NulandSukochev2023,vNvS21a}. Under reasonable assumptions, these identities are
    \begin{align}
        f(H+V)-f(H)&=T_{f^{[1]}}^{H+V,H}(V),\label{eq:lowneridentity}\\
        [f(H),a]&=T_{f^{[1]}}^{H,H}([H,a])\label{eq:commutator},
    \end{align}
    and their higher-order analogues, cf. Proposition~\ref{P:UMOIcom}.

    Applying an identity like~\eqref{eq:commutator} to the JLO cocycle where $H= D^2$ gives an unbounded argument $[D^2, a]$ in the multiple operator integral, for which the current theory of multiple operator integrals~\cite{AleksandrovPeller2022,ACDS,FrankPushnitski2019, Peller2016, SkripkaTomskova2019} is not applicable. In literature on noncommutative geometry, ad-hoc arguments have been made to handle this~\cite{CPRS1, ConnesMoscovici1995, Higson2003}, but never by generalising multiple operator integrals techniques to deal with unbounded arguments. 

\subsection{Summary of Main Results}
    The pseudodifferential calculus we use in our approach can be regarded as the operator-theoretical skeleton behind any pseudodifferential calculus. It is therefore difficult to pin down its historical origins.  The abstract Sobolev spaces in this theory, also called Hilbert scales, already appeared in their current form in work by S. G. Kre\u{\i}n, see the reviews~\cite{Mitjagin1961, KreinPetunin1966} and references therein. The operators themselves seem to have been pioneered by Connes--Moscovici, Higson and Guillemin~\cite{ConnesMoscovici1995, Guillemin1985,Higson2003}, and developed further by Br\"uning and Lesch~\cite{BruningLesch2001}. A succinct overview of the properties of this calculus is given in~\cite{Uuye2011}.

\begin{defn}\label{def:pseudocalc}
    Let $\Theta$ be a possibly unbounded invertible positive self-adjoint operator on a separable Hilbert space $\H$. Define the Hilbert spaces $\H^s := \overline{\dom \Theta^s}^{\| \cdot \|_s}$ for $s \in \mathbb{R}$ where $\|\phi\|_{s} : = \| \Theta^s \phi\|$ -- though taking this closure is not necessary for $s\geq 0$. We write $\H^\infty := \bigcap_{s \geq 0} \H^s$, which is dense in $\H$. We say that a linear operator $A: \H^\infty \to \H^\infty$ is in the class $\op^r(\Theta)$ (it has analytic order $\leq r$) if $A$ extends to a continuous operator
    \[
    \overline{A}^{s+r,r}: \H^{s+r} \to \H^{s}
    \]
    for all $s\in \mathbb{R}$. If no confusion can arise, we often write
    \[
    A: \H^{s+r} \to \H^s,
    \]
    $\op(\Theta) := \bigcup_{r\in \mathbb{R}} \op^r(\Theta)$, and $\op^{-\infty}(\Theta):= \cap_{r\in \mathbb{R}}\op^r(\Theta)$. Finally, we define $\OP^r(\Theta)\subseteq \op^r(\Theta)$ as those $A \in \op^r(\Theta)$ for which $\delta^n_\Theta(A) \in \op^r(\Theta)$ for each $n \geq 0$, where $\delta_\Theta(A) := [\Theta, A]$. 
\end{defn}

Both $\op(\Theta)$ and $\OP(\Theta):= \bigcup_{r \in \R} \OP^r(\Theta)$ form a filtered algebra, as $\op^r(\Theta) \cdot \op^t(\Theta) \subseteq \op^{r+t}(\Theta)$ and $\OP^r(\Theta)\cdot \OP^t(\Theta) \subseteq \OP^{r+t}(\Theta)$. 

Commonly, this paradigm is also used when dealing with unbounded operators
\[
T: \dom(T) \to \H.
\]
In this case, one writes $T \in \op^r(\Theta)$ if $\H^\infty \subseteq \dom(T)$, $T(\H^\infty) \subseteq \H^\infty$, and
\[
T\big|_{\H^\infty} \in \op^r(\Theta).
\]
Conversely, for $T \in \op^r(\Theta)$ with $r>0$,
\[
\overline{T}^{r,0}: \H^r \subseteq \H \to \H
\]
can be interpreted as an unbounded operator.
Furthermore, note that $\op^{-r}(\Theta) \subseteq \mB(\H)$ for $r \geq 0$ in the sense that for $A \in \op^{-r}(\Theta)$ we have $\overline{A}^{-r,0}|_{\H} \in B(\H)$. Similarly if $\Theta^{-1}\in \mathcal{L}_s(\H)$, a Schatten--von Neumann class, then $\op^{-s}(\Theta) \subseteq \mathcal{L}_{1}(\H)$.

Let us provide a few examples of this abstract pseudodifferential calculus.
\begin{itemize}
    \item Taking $\Theta = (1-\Delta)^{1/2}$ on $L_2(\mathbb{R}^d)$, where $\Delta=\sum_{j=1}^d \partial_j^2$ is the Laplace operator, gives the classical (Bessel potential) Sobolev spaces $\H^s = W_2^s(\mathbb{R}^d)$. (Pseudo)differential operators of order $k$, and (unbounded) Fourier multipliers $T_\phi$ with symbols $|\phi(\xi)|\lesssim (1+|\xi|)^k,\, \xi\in \mathbb{R}^d$ are contained in $\OP^k((1-\Delta)^{1/2})$. Note though that $\OP((1-\Delta)^{1/2})$ is a larger class than this, which for example also contains translation operators.
    \item  For spectral triples $(\A, \H, D)$ taking $\Theta= (1+D^2)^{1/2}$ recovers the calculus of Connes and Moscovici as used in noncommutative geometry.
    \item Taking $\Theta = (1-\Delta)^{1/2}$, where $\Delta$ is the sub-Laplacian on a stratified Lie group, gives the Sobolev spaces defined by Folland and Stein~\cite{FollandStein1974,RothschildStein1976}. 
    \item Related to the previous example is the anharmonic oscillator $\Theta^2 = 1- \Delta^{2l} + |x|^{2k} $ on $\R^d$ for integers $l,k\geq 1$ and generalisations thereof, which define Sobolev spaces and a pseudodifferential calculus that appear in the study of sub-Laplacian operators on stratified Lie groups too~\cite{ChatzakouDelgado2021}. The special case where $\Theta^2$ is the harmonic oscillator gives Shubin's Sobolev spaces $Q^s(\R^d)$~\cite[Section~IV.25]{Shubin2001} and the $\Gamma$-pseudodifferential operators in~\cite[Chapter~2]{NicolaRodino2010}, see also~\cite{BongioanniTorrea2006}.
    \item In recent work by Androulidakis, Mohsen and Yuncken, Sobolev spaces are constructed with a similar procedure in order to prove the Helffer--Nourrigat conjecture~\cite{AndroulidakisMohsen2022}. 
    \item A similar calculus has been constructed for quantum Euclidean spaces~\cite{GaoJunge2022}.
    \item Finally, we note the case where $\Theta = 1_{\H}$, which gives that $\H^\infty = \H^s = \H$ for all $s \in \R$, and $\op^r(1_{\H}) = B(\H)$, $r \in \R$.
\end{itemize}

In this formalism, we prove a generalisation of the MOI framework by Peller and Azamov--Carey--Dodds--Sukochev~\cite{ACDS, Peller2006}. Recall the following definition, cf.~\cite[Chapter~4]{Schmudgen2012},~\cite[Section~3.2]{SkripkaTomskova2019}. Throughout this paper, we use the notation $\langle x \rangle := (1+|x|^2)^{1/2}$. 

\begin{defn}\label{def:LinftySpectral}
Let $E$ be a spectral measure on $\R$ with the Borel sigma algebra. For a Borel measurable function $f:\R\to\C$ we define the essential supremum seminorm
$$\|f\|_{L_\infty(E)}:=\sup\{y\in\R~:~E(|f|^{-1}((y,\infty)))=0\},$$
which defines $L_\infty(E)$ in the usual way, namely as the quotient of the set of measurable functions with finite seminorm, by the set of those of zero seminorm. In the same way we define $L_\infty^\beta(E)$ by the seminorm
$$\|f\|_{L_\infty^\beta(E)}:=\|x\mapsto f(x)\langle x\rangle^{-\beta}\|_{L_\infty(E)},$$
for any $\beta\in\R$, where $\langle x \rangle = (1+x^2)^{1/2}$.
\end{defn}

\begin{thm}\label{T:MainMOIConstruction}
Let $H_0,\ldots H_n$ be self-adjoint operators on $\H$ with spectral measures $E_j$, and let $\phi:\R^{n+1}\to\C$ be of the form
\[
\phi(\lambda_0,\ldots,\lambda_n)=\int_\Omega a_0(\lambda_0,\omega)\cdots a_n(\lambda_n,\omega) d\nu(\omega),
\]
for a finite measure space $(\Omega,\nu)$ and measurable functions $a_j:\R\times\Omega\to\C$ such that $(x,\omega)\mapsto a_j(x,\omega)\langle x \rangle^{-\beta_j}$ is $E_j \times \nu$-a.e. bounded for $\beta_j \in \R$. Suppose that we have $a_j(H_j,\omega)\in\op^{k_j}(\Theta)$, $k_j \in \R$, and 
\begin{equation}\label{eq:normboundfunctcalc}
    \|a_j(H_j,\omega)\|_{\H \indices{^{s+k_j}}\to\H\indices{^s}}\leq C_{s,H_j}\|a_j(\cdot,\omega)\|_{L^{\beta_j}_{\infty}(E_j)}
\end{equation}
for every $j \in \{0, \ldots, n\}$, $s\in \R$, and $\omega\in\Omega$, and certain constants $C_{s,H_j}\in\R$. Then the integral
\[
T_\phi^{H_0,\ldots,H_n}(X_1,\ldots,X_n)\psi:=\int_\Omega a_0(H_0,\omega)X_1 a_1(H_1,\omega)\cdots X_n a_n(H_n,\omega)\psi\,d\nu(\omega), \quad \psi\in\H^\infty,
\]
for $X_1,\ldots,X_n\in\op(\Theta)$, converges as a Bochner integral in $\H^s$ for every $s \in \R$, and defines an $n$-multilinear map $T_\phi^{H_0,\ldots,H_n}:\op^{r_1}(\Theta)\times\cdots\times\op^{r_n}(\Theta)\to\op^{\sum_j r_j + \sum_j k_j}(\Theta)$ depending on $\Omega$ and $a_0,\ldots,a_n$ only through the symbol~$\phi$. Specifically, for $s\in \R$ we have the estimate
   \begin{align*}
        \big \|T_\phi^{H_0,\ldots, H_n}(X_1, \ldots,X_n)\big\|_{\H^{s+\sum_{j}r_j + \sum_j k_j} \to \H^{s}} \lesssim \prod_{j=1}^n \|X_j\|_{\H^{s_{j}+r_j} \to \H^{s_{j}}} \int_\Omega \prod_{j=0}^n \| a_j(\cdot, \omega)\|_{L^{\beta_j}_\infty(E_j)} d|\nu|(\omega),
    \end{align*}
    for some $s_1, \ldots, s_n \in \R$.
\end{thm}

For $\Theta = 1_{\H}$ and $\beta_0 = \cdots = \beta_n = 0$, it is immediate that $a_j(H_j, \omega) \in \op(1_{\H})= B(\H)$ and
\[
\|a_j(H_j, \omega)\|_{\H^{s} \to \H^s} = \|a_j(\cdot, \omega)\|_{L^0_\infty(E_j)},
\]
and the above theorem reduces to Peller's and Azamov--Carey--Dodds--Sukochev's construction of MOIs~\cite{ACDS, Peller2006, Peller2016}. For general $\Theta$, the question is for which self-adjoint operators $H$ and functions $f$ we have $f(H) \in \op(\Theta)$ with the required norm estimate~\eqref{eq:normboundfunctcalc}. In particular, for $H \in \op^h(\Theta)$, $h>0$, the extension $\overline{H}^{r,0}$ can be considered an unbounded operator on $\H$ and hence (if it is self-adjoint) we might ask whether $f(\overline{H}^{r,0}) \in \op(\Theta)$.
We draw inspiration from the traditional theory of pseudodifferential operators and consider elliptic symmetric operators in $\op(\Theta)$.

We say that an operator $A \in \op^r(\Theta)$, $r\geq0$, is \textit{symmetric} or \textit{self-adjoint with domain $\H^r$} respectively, if the extension
\[
\overline{A}^{r,0}: \H^r \subseteq \H^0 \to \H^0,
\]
regarded as an unbounded operator on $\H^0$, is a symmetric or self-adjoint operator. See Section~\ref{S:FunctCalc} and Appendix~\ref{S:AppEllipticAdj} for an investigation of these properties.

We say that $A\in \op^r(\Theta)$ is \textit{$\Theta$-elliptic} if there exists a parametrix $P \in \op^{-r}(\Theta)$ such that
\begin{align*}
    AP &= 1_{\H^\infty} + R_1;\\
    PA &= 1_{\H^\infty} + R_2,
\end{align*}
where $R_1, R_2 \in \op^{-\infty}(\Theta).$ This is similar to \cite[Definition 2.1]{Guillemin1985}. The notion of $\Theta$-ellipticity depends on $\Theta$ and on the order $r \in \R$. We simply abbreviate this to `$A \in \op^r(\Theta)$ is $\Theta$-elliptic'. Note that when $\Theta = X_1^2+\cdots+X_n^2$ is the H\"ormander Laplacian~\cite{Hormander1967}, $\Theta$-ellipticity coincides precisely with what is usually called `maximal hypoellipticity' in this context, cf.~\cite{AndroulidakisMohsen2022}. 

Properties of $\Theta$-elliptic operators are discussed in more depth in Section~\ref{S:FunctCalc}, where it is shown that symmetric $\Theta$-elliptic operators in $\op^r(\Theta)$, $r\geq 0$ are self-adjoint with domain $\H^r$. Therefore, when we mention the spectrum or spectral measure of such operators, we are referring to the spectrum or spectral measure of the extension $\overline{A}^{r,0}$.

\begin{thm}\label{T:MainFunctCalc}
   Let $A\in \op^r(\Theta)$, $r >0$, be $\Theta$-elliptic and symmetric, and let $E$ denote its spectral measure. If $f \in L^\beta_\infty(E), \beta \in \R,$ then 
        \[
        f(\overline{A}^{r,0}) \in \op^{\beta r}(\Theta),
        \]
        and we simply write $f(A) := f(\overline{A}^{r,0})$.
        More precisely,
        \[
        \|f(A)\|_{\H^{s+\beta r} \to \H^s} \leq C_{s,A} \|f \|_{L^\beta_\infty(E)}.
        \]
\end{thm}

This theorem is proved in Section~\ref{S:FunctCalc}. The result might be surprising: for traditional elliptic pseudodifferential operators, a similar functional calculus is only constructed with \textit{smooth} functions $f$ in $S^\beta(\R)$~\cite[Theorem 2.4]{Bony2013}, \cite{Robert1982,Strichartz1972}, \cite[Chapter~XII]{Taylor1981}. A functional calculus for $\op^0(\Theta)$ is developed in Appendix~\ref{S:AppFunctCalc0}.

\begin{defn}\label{def:SbetaTbeta} For $I\subseteq \R$ an interval and $\beta \in \R$, we define $S^\beta(I)$ as the class of smooth functions $f: I \to \C$ such that
\[
\|f\|_{S^\beta(I), k} := \sup_{x\in I} |f^{(k)}(x)| \langle x \rangle^{k-\beta} < \infty, \quad k \in \Z_{\geq 0}.
\]
The quantities above are seminorms. Similarly, we define $T^{\beta}(I)$ as the space of smooth functions $f: I \to \C$ such that 
    \[
    \|f\|_{T^\beta(I), k}:= \int_I |f^{(k)}(x)|\langle x\rangle^{k-\beta-1}dx < \infty, \quad k \in \Z_{\geq 0}.
    \]
\end{defn}
    In case $I = \R$, we note the inclusions
    \[
    \bigcup_{\alpha<\beta} S^{\alpha}(\R) \subsetneq T^{\beta}(\R) \subsetneq S^\beta(\R) \subsetneq L^\beta_\infty(E)
    \]
    for any spectral measure $E$.

As a combination of Theorem~\ref{T:MainMOIConstruction} and Theorem~\ref{T:MainFunctCalc}, we now obtain the following.
\begin{cor}\label{C:MOOIforNCG1}
    Let $H_0, \ldots, H_n$ be such that each $H_i \in \op^{h_i}(\Theta), h_i > 0,$ is symmetric and $\Theta$-elliptic with spectral measure $E_i$.
    Let $\phi: \mathbb{R}^{n+1} \to \mathbb{C}$ be such that
    \[
    \phi(\lambda_0, \ldots, \lambda_n) = \int_\Omega a_0(\lambda_0, \omega) \cdots a_n(\lambda_n, \omega) d\nu(\omega),
    \]
    for a finite measure space $\Omega$ and measurable functions $a_i: \R \times \Omega \to \C$ such that $(x,\omega)\mapsto a_i(x,\omega)\langle x \rangle^{-\beta_i}$ is $E_i \times \nu$-a.e. bounded for $\beta_i \in \R$.
    Then Theorem~\ref{T:MainMOIConstruction} applies, and for operators $X_i \in \op^{r_i}(\Theta)$, $r:= \sum_{i=1}^n r_i$, we have that 
    \[
    T_{\phi}^{H_0, \ldots, H_n}(X_1, \ldots, X_n) \in \op^{r+\sum_{i=0}^n \beta_i h_i}(\Theta)
    \]
    independent of the chosen representation of $\phi$.
\end{cor}

In order to make use of this construction, we need to ensure that divided differences $f^{[n]}$ under appropriate conditions on $f$, admit such integral decompositions
\[
f^{[n]}(\lambda_0, \ldots, \lambda_n) = \int_\Omega a_0(\lambda_0, \omega) \cdots a_n(\lambda_n, \omega) d\nu(\omega),
\]
where each $a_i(x,\omega)\langle x\rangle^{-\beta_i} $ is $E_i \times \nu$-a.e. bounded for $\beta_i \in \R$. This is the main achievement of Section~\ref{S:FunctSpac}, where we employ the theory of almost analytic extensions to obtain integral decompositions of divided differences in a novel way. In summary, we prove that if $f \in C^{n+2}(\R)$ and $\|f\|_{T^\alpha(\R),k}<\infty$ for $k=0, \ldots, n+2$, then the divided difference $f^{[n]}$ can be decomposed in a finite sum of functions $\phi$ which themselves admit an integral decomposition
\[
\phi(\lambda_0, \ldots, \lambda_n) = \int_\Omega a_0(\lambda_0, \omega) \cdots a_n(\lambda_n, \omega) d\nu(\omega),
\]
with each $a_i(\cdot, \omega) \in L^{\beta_i}_\infty(\R)$. Depending on whether $\alpha$ is in the regime $(-\infty, -1]$, $[-1,n]$ or $[n,\infty)$, we have varying control over the coefficients $\beta_i \in \R$, which sum up to $\alpha-n$. The precise statement is given in Lemma~\ref{L:DivDifT}. Corollary~\ref{C:MOOIforNCG1} and Lemma~\ref{L:DivDifT} can be assembled into the following theorem, which can be regarded as the main theorem of this paper.

     \begin{thm}\label{T:MainThmForReal}
     Let $H_0, \ldots, H_n$ be such that each $H_i \in \op^{h_i}(\Theta), h_i > 0,$ is symmetric and $\Theta$-elliptic with spectral measure $E_i$. Let $f \in C^{n+2}(\R)$ such that $\|f\|_{T^\beta(\R), k} <\infty$, $k=0, \ldots, n+2$ for some $\beta \in \R$. There exists $C \in \R$ depending only on $h_0,\ldots,h_n,\beta,n$, such that 
     \[
     T^{H_0, \ldots, H_n}_{f^{[n]}}: \op^{r_1}(\Theta) \times \cdots \times \op^{r_n}(\Theta) \to  \op^{C+\sum_{j=1}^n r_j}(\Theta),
     \]
     is a well-defined multilinear map. For $X_i \in \op^{r_i}(\Theta)$, $r:= \sum_{i=1}^n r_i$, we have that there exist  $s_1, \ldots, s_n \in \R$ such that
    \[
      \| T^{H_0, \ldots, H_n}_{f^{[n]}}(X_1, \ldots, X_n)\|_{\H^{s+r+ C} \to \H^s} \leq C_{s,H_0, \ldots, H_n} \left(\sum_{k=0}^{n+2} \|f\|_{T^\beta(\R), k} \right) \prod_{i=1}^n \|X_i\|_{\H^{s_i+r_i} \to \H^{s_i}}.
     \]

     In the case where $h_0 = \cdots = h_n=:h$, the constant $C$ takes the value $(\beta-n)h$, so that
     \[
     T^{H_0, \ldots, H_n}_{f^{[n]}}(X_1, \ldots, X_n) \in \op^{r+(\beta-n)h}(\Theta).
     \]
\end{thm}

\begin{rem}
    The functional calculus for $\Theta$-elliptic symmetric operators in $\op(\Theta)$ in Theorem~\ref{T:MainFunctCalc} can be extended to self-adjoint operators $H$ (not necessarily in $\op(\Theta)$) that strongly commute with some $\Theta$-elliptic operator in $\op(\Theta)$ of positive order. See Proposition~\ref{P:BigFunctCalc} for details. The multiple operator integral theory presented in Theorem~\ref{T:MainThmForReal} therefore also extends to the case where the operators $H_i$ are of this type.
\end{rem}

Theorem~\ref{T:MainThmForReal} can of course be applied when all operators $H_i$ and $X_i$ are of the class $\OP(\Theta) \subseteq \op(\Theta)$, and it may be of interest to know when the resulting operator is again of the class $\OP(\Theta)$. In case the symbol is a divided difference we have the following result, proved in Section~\ref{SS:Identities}.

\begin{thm}\label{T:MOIpsi}
    If $H_0, \ldots, H_n \in \OP^{h}(\Theta)$ are symmetric and $\Theta$-elliptic of the same order $h>0$ and $X_j \in \OP^{r_j}(\Theta)$, $r:=\sum_{j=1}^n r_j$, and if $f \in T^{\beta}(\R),$ then the multiple operator integral $T^{H_0, \ldots, H_n}_{f^{[n]}}(X_1, \ldots, X_n) \in \op^{(\beta - n)h + r}(\Theta)$ whose existence follows from Theorem~\ref{T:MainThmForReal}, is again an element of $\OP(\Theta)$, i.e.
    \[
    T^{H_0, \ldots, H_n}_{f^{[n]}}(X_1, \ldots, X_n) \in \OP^{(\beta - n)h + r }(\Theta).
    \]
    In particular, the $n=0$ case gives that
    \[
    f(H_0) \in \OP^{\beta h}(\Theta).
    \]
\end{thm}

\begin{rem}
    The operators $T_\phi^{H_0,\ldots, H_n}(X_1, \ldots, X_n)$ constructed in Theorem~\ref{T:MainMOIConstruction}, Corollary~\ref{C:MOOIforNCG1} and Theorem~\ref{T:MOIpsi} depend on the symbol $\phi$ only through the function value of $\phi$ on the spectra $\sigma(H_0) \times \cdots \times \sigma(H_n) \subseteq \R^n$. Hence all statements regarding the operators of the form $T^{H}_{f^{[n]}}(X_1, \ldots, X_n)$ where $f \in C^{n+2}(\R)$ or $f \in T^\beta(\R)$ remain valid if $f$ is in $C^{n+2}(I)$ or $T^\beta(I)$ respectively, where $I\subseteq \R$ is an open neighbourhood of the spectrum of $H$.
\end{rem}

For a \textit{regular} spectral triple, which means that $\A, [D,\A] \subseteq \OP^{0}((1+D^2)^{1/2})$, the applicability of the above theorems is obvious.

Given these constructions for multiple operator integrals, various identities follow quite immediately, as will be shown in Section~\ref{S:MOIIdentsAndAsympt}. Of particular interest are the following expansions. We say that~\cite{Higson2003}
\[
A \sim \sum_{k=0}^\infty A_k
\]
for $A, A_k \in \op(\Theta)$, if we have
\begin{equation}\label{eq:pdoexpansion}
A - \sum_{k=0}^N A_k \in \op^{m_N}   (\Theta) 
\end{equation}
with $m_N \downarrow -\infty$.

\begin{thm}\label{T:AsypExpansions}
    \begin{enumerate}
        \item Let $f \in T^\beta(\R)$, take $H \in \op^h(\Theta)$, $h >0$ $\Theta$-elliptic and symmetric, let $V \in \op^r(\Theta)$ be symmetric. If the order of the perturbation $V$ is strictly smaller than that of $H$, i.e. $r < h$, we have
    \begin{align*}
    f(H + V) \sim &\sum_{n=0}^\infty T^{H}_{f^{[n]}}(V, \ldots, V).
    \end{align*}
    \item Let $f \in T^\beta(\R)$, take $H \in \op^h(\Theta)$, $h >0$ $\Theta$-elliptic and symmetric, and take $X_i \in \op^{r_i}(\Theta)$. If there exists some $\varepsilon > 0$ such that $\delta^n_{H}(X_i)\in \op^{r_i + n (h-\varepsilon)}(\Theta)$ for all $i$ and all $n\geq 0$, then we have an asymptotic expansion
    \begin{align*}
          T^{H}_{f^{[n]}}(X_1, \dots, X_n) \sim \sum_{m=0}^{\infty} \sum_{m_1 + \dots + m_n = m} \frac{C_{m_1, \dots, m_n}}{(n+m)!} \delta_H^{m_1}(X_1) \cdots \delta_H^{m_n}(X_n) f^{(n+m)}(H) ,
    \end{align*}
    where 
    \[
    C_{m_1, \dots, m_n} :=  \prod_{j=1}^{n} \binom{j+m_1 + \dots + m_{j}-1}{m_j}.
    \]
    \item Under the assumptions of the first two parts combined, we obtain
    \begin{align*}
          f(H+V) \sim  \sum_{n,m=0}^\infty \sum_{m_1 + \dots + m_n = m} \frac{C_{m_1, \dots, m_n}}{(n+m)!} \delta_H^{m_1}(V) \cdots \delta_H^{m_n}(V) f^{(n+m)}(H).
    \end{align*}
    \end{enumerate}
\end{thm}

This theorem is proved through Theorem~\ref{T:Taylor}, Proposition~\ref{P:Expansion} and Corollary~\ref{C:CombiExpansion}. The first part of the theorem above should be interpreted as a Taylor expansion. Special cases of the second part of the theorem have proved essential in proofs of the local index formula~\cite{ConnesMoscovici1995, CPRS1}, where they have been derived in a case by case basis. A version for $C^*$-algebras has also appeared in~\cite{lesch2017divided}, and a version for classical pseudodifferential operators in~\cite{Paycha2007}. The condition that $\delta^n_H(X_i) \in \op^{r_i + n(h-\varepsilon)}(\Theta)$ is for example satisfied for $H = \Theta$ and $X_i \in \OP(\Theta)$, with $\varepsilon = 1$.

Apart from its use in finding asymptotic expansions of trace formulas, a seemingly disparate application of the above theorem is given by the following corollary which is also used throughout the literature on noncommutative geometry. In particular, a version appears as Theorem B.1 in Connes--Moscovici~\cite{ConnesMoscovici1995}. It has been used repeatedly in the noncommutative geometry literature, see for example~\cite{CPRS2, Higson2003}, and for a more recent example~\cite{Rodsphon2015}.

\begin{cor}
    For $X \in \OP^r(\Theta)$ and $f \in T^\beta(\R_{>\varepsilon})$ we have
    \[
    [f(\Theta), X] \sim \sum_{k=1}^\infty \frac{1}{k!}\delta_\Theta^k(X) f^{(k)}(\Theta).
    \]
    In particular,
    \[
    [\Theta^\alpha, X] \sim \sum_{k=1}^\infty \binom{\alpha}{k}\delta_\Theta^k(X) \Theta^{\alpha - k}, \quad \alpha \in \C,
    \]
    and
    \[
    [\log(\Theta), X] \sim \sum_{k=1}^\infty \frac{ (-1)^{k-1}}{k} \delta_\Theta^k(X) \Theta^{- k},
    \]
    and we have that $[\Theta^\alpha, X] \in \OP^{r+\Re(\alpha)-1}(\Theta)$ and $[\log(\Theta),X] \in \OP^{r-1}(\Theta)$.
\end{cor}
\begin{proof}
We first reduce to $r=0$ by writing $\tilde{X} = X\Theta^{-r} \in \OP^0(\Theta)$ and noting that
\[
[f(\Theta), X] = [f(\Theta), \tilde{X}]\Theta^r, \quad \delta_\Theta^k(X) = \delta^k_\Theta(\tilde{X})\Theta^r.
\]
For $r=0$, the assertion is a simple combination of the second part of Theorem~\ref{T:AsypExpansions} with the identity
\[
[f(\Theta), X] = T^{\Theta}_{f^{[1]}}([\Theta, X])
\]
mentioned earlier, which will be proved in Proposition~\ref{P:UMOIcom}. 
\end{proof}

Finally, the expansions in Theorem~\ref{T:AsypExpansions} can be refined to provide trace expansions. We apply this to prove a pertubative expansion of the spectral action in noncommutative geometry~\cite[Chapter~7]{Suijlekom2015}\cite{EcksteinIochum2018}, and a result on the existence of asymptotic expansions of the heat trace of the perturbed Dirac operator.
For the purpose of understanding the following theorem, one only needs to know that a regular $s$-summable spectral triple $(\A, \H, D)$ consists of a Hilbert space $\H$, an algebra of bounded operators $\A$ and an unbounded self-adjoint operator $D$. The triple being regular (sometimes called smooth or $QC^\infty$) means that if we write $\Theta = (1+D^2)^{1/2}$, we have
\[
\A, [D,\A] \subseteq \OP^0(\Theta),
\]
and the $s$-summability means that $\Theta^{-1} \in \mathcal{L}_s$. For a more indepth overview, see~\cite{CGRS2, CareyPhillips2011, EcksteinIochum2018, Higson2003, Rennie2003, Suijlekom2011}.

\begin{thm}\label{T:AsympExpNCG}
Let $(\A,\H,D)$ be a regular $s$-summable spectral triple, $s> 0$. Let $V, P \in \mB$, $V$ self-adjoint and bounded, where $\mB$ is the algebra generated by $\A$ and $D$. 
   Then for all $f \in T^{\beta}(\R)$ with $\beta < -s$, the expressions
\begin{align*}
\Tr(f(tD + tV)), \qquad \Tr(Pe^{-t(D+V)^2}) \quad \text{and} \quad \Tr(Pe^{-t|D+V|})
\end{align*}
admit an asymptotic expansion as $t \downarrow 0$ given respectively by
\[
    \Tr(f(tD+tV)) = \sum_{n=0}^{N} \sum_{m=0}^{N}   \sum_{m_1 + \dots + m_n = m} t^{n+m} \frac{C_{m_1, \dots, m_n}}{(n+m)!} \Tr\big( \delta_D^{m_1}(V) \cdots \delta_D^{m_n}(V) f^{(n+m)}(tD)\big) + O(t^{N+1-s}),
\]
where $C_{m_1, \dots, m_n}$ is the same as in Theorem~\ref{T:AsypExpansions}, 
\[
\Tr(Pe^{-t(D+V)^2}) = \sum_{n = 0}^{N} \sum_{m=0}^{N} \sum_{m_1 + \dots + m_n = m}  (-t)^{n+m} \frac{C_{m_1, \dots, m_n}}{(n+m)!} \Tr(PA^{(m_1)}\cdots A^{(m_n)}\exp(-tD^2)) + O(t^{\frac{N+1-s}{2}}),
\]
where $A := DV + VD + V^2$, and $A^{(m)} := \delta^m_{D^2}(A)$, 
and
\begin{align*}
    &\Tr(Pe^{-t|D+V|}) \\
    &= \sum_{n = 0}^{N} \sum_{m=0}^{N} \sum_{m_1 + \dots + m_n = m}  (-t)^{n+m} \frac{C_{m_1, \dots, m_n}}{(n+m)!} \Tr(P\delta_{|D|}^{m_1}(B)\cdots \delta_{|D|}^{m_n}(B)\exp(-t|D|)) + O(t^{(N+1)(1-\varepsilon) - s}),
\end{align*}
where $B := |D+V| - |D|$ and $\varepsilon > 0$ can be picked arbitrarily small.
\end{thm}

This theorem is proved in Section~\ref{SS:AsympExp}. In the context the spectral action, this expansion is closely related to~\cite{Skripka2014, Suijlekom2011}. Both these cited papers and Theorem~\ref{T:AsympExpNCG} are based on a Taylor expansion of the spectral action, but the expansion in Theorem~\ref{T:AsympExpNCG} goes one step further by applying the second part of Theorem~\ref{T:AsypExpansions}. Observe that we do not assume that $f$ is the Laplace transform of a measure on $\R_+$, which is a common condition in the literature.

In Section~\ref{S:Basics} we cover some preliminary material for this paper and prove Theorem~\ref{T:MainMOIConstruction}. Section~\ref{S:FunctCalc} proves the existence of a functional calculus for $\Theta$-elliptic symmetric pseudodifferential operators, and Section~\ref{S:FunctSpac} shows that divided differences of functions in $T^\beta(\R)$ satisfy the conditions required for symbols by Theorem~\ref{T:MainMOIConstruction}, thus completing the proof of Theorem~\ref{T:MainThmForReal}. Finally, in Section~\ref{S:MOIIdentsAndAsympt} we prove some identities for our construction of a MOI and we finish with a number of results on asymptotic expansions. Appendix~\ref{S:AppEllipticAdj} contains a brief study on adjoints of abstract pseudodifferential operators, Appendix~\ref{S:AppFunctCalc0} gives a functional calculus for $\op^0(\Theta)$ and Appendix~\ref{S:AppCombinComp} contains the (combinatorial) proof of one of the expansions in Section~\ref{S:MOIIdentsAndAsympt}.

\textbf{Acknowledgements}

\noindent The authors wish to thank Fedor Sukochev and Dmitriy Zanin for helpful discussions and comments. We thank Bruno Iochum for communicating the question which inspired part of this paper, and for correspondence which improved the manuscript. Dmitriy Zanin provided help with the proof of Proposition~\ref{P:Invertop}.

\numberwithin{thm}{section}

\section{Multiple Operator Integrals as Pseudodifferential Operators}\label{S:Basics}
In this section we prove Theorem~\ref{T:MainMOIConstruction}. The proof is a subtle modification of the proof presented in~\cite{ACDS, Peller2006}, so the material presented here is mostly an adaptation of known results.

\subsection{Operator Integrals}
We start with standard definitions and results on measurability and integrability of operator valued functions. Throughout the paper we fix a separable Hilbert space $\H$ and a positive invertible operator $\Theta$ on $\H$, which yields separable Hilbert spaces $\H^s$ by Definition~\ref{def:pseudocalc}.
\begin{defn}
    Let $\H_0, \H_1$ be separable Hilbert spaces and let $(\Omega, \Sigma, \nu)$ be a measure space with complex measure. A function $f: \Omega \to B(\H_1, \H_0)$ is called weak operator topology measurable (weakly measurable for short) if for all $\eta \in \H_0$, $\xi \in \H_1$ the scalar-valued function
    \[
    \omega \mapsto \langle \eta, f(\omega) \xi \rangle_{\H_0}, \quad \omega \in \Omega,
    \]
    is measurable. Similarly, $f$ is said to be weak operator topology integrable if for all $\xi$ and $\eta$ the above map is integrable.
\end{defn}

\begin{lem}{\cite[Lemma~1.4.2]{LMSZVol2}}\label{L:PettisIntegral}
Let $\H_0, \H_1, (\Omega, \Sigma, \nu)$ be as above, and let $f: \Omega \to B(\H_1, \H_0)$ be weakly measurable. Then the norm function
    \[
    \omega \mapsto \|f(\omega)\|_{\H_1 \to \H_0}, \quad \omega \in \Omega,
    \]
    is measurable. If moreover
    \[
    \int_\Omega \| f(\omega)\|_{\H_1 \to \H_0} d|\nu|(\omega) < \infty,
    \]
    then there exists a unique $I_f \in  B(\H_1, \H_0)$ such that
    \[
    \langle \eta, I_f \xi \rangle_{\H_0} = \int_\Omega \langle \eta, f(\omega)\xi \rangle_{\H_0} d\nu(\omega), \quad \eta \in \H_0, \xi \in \H_1,
    \]
    and
    \[
    \| I_f \|_{\H_1 \to \H_0} \leq \int_\Omega \| f(\omega)\|_{\H_1 \to \H_0} d|\nu|(\omega).
    \]
    We then write $I_f = \int_\Omega f(\omega) d\nu(\omega)$.
\end{lem}

\begin{prop}{\cite[Lemma~3.11]{DoddsDodds2020II}}\label{P:WeakMeasFunctCalc}
    Let $(\Omega, \Sigma, \nu)$ be a $\sigma$-finite measure space, let $a:\R \times \Omega \to \C$ be measurable and bounded, and let $H$ be an (unbounded) self-adjoint operator on $\H$. Then
    \[
    \omega \mapsto a(H,\omega)
    \]
    is weakly measurable.
\end{prop}
\begin{proof}
Though~\cite[Lemma~3.11]{DoddsDodds2020II} is only formulated for bounded $H$, the unbounded case follows with the same proof. 
\end{proof}

\begin{lem}\label{L:WeaklyMeasComp}
    Let $\H_0, \ldots, \H_{2n+1}$ be separable Hilbert spaces, $X_i \in B(\H_{2i}, \H_{2i-1})$, and let $f_i:\Omega \to B(\H_{2i+1}, \H_{2i} )$ be weakly measurable functions. Then
    \begin{align*}
        \Omega & \to B(\H_{2n+1}, \H_0)\\
        \omega &\mapsto f_0(\omega) X_1 f_1(\omega) \cdots X_n f_n(\omega)
    \end{align*}
    is weakly measurable. Furthermore, if 
    \[
    \int_\Omega \| f_0(\omega)\|_{\H_1 \to \H_0} \cdots \|f_n(\omega)\|_{\H_{2n+1}\to \H_{2n}} d|\nu|(\omega) < \infty,
    \]
    the map
    \begin{align*}
    B(\H_{1}, \H_{0}) \times \cdots \times B(\H_{2n-1}, \H_{2n-2}) &\to B(\H_{2n}, \H_0)\\
    (X_1, \ldots, X_n) \quad &\mapsto \int_\Omega f_0(\omega) X_1f_1(\omega) \cdots X_n f_n(\omega) d\nu(\omega),
\end{align*}
whose existence follows from Lemma~\ref{L:PettisIntegral}, is $SOT$-continuous when restricted to the unit ball in each argument $B(\H_{2i}, \H_{2i-1})$.
\end{lem}
\begin{proof}
    The first part of the lemma is a consequence of the fact that the pointwise product of weakly measurable functions is weakly measurable, see~\cite[Lemma~3.7]{DoddsDodds2020II}.
    
    The $SOT$-continuity follows from the joint continuity of the multiplication
    \[
     (X_1, \ldots, X_n) \mapsto  a_0(H_0, \omega) X_1 a_1(H_1, \omega) \cdots X_n a_n(H_n, \omega)
    \]
    in the strong operator topology when restricting to the unit balls~\cite{Blackadar2006}, in combination with the Dominated Convergence Theorem for the Bochner integral of Hilbert space-valued functions~\cite[Corollary~III.6.16]{DunfordSchwartzI}.
\end{proof}

We use the notation $\| \cdot \|_\infty$ for the essential supremum.

\begin{lem}\label{L:EssSupMeasurable}
    Let $(\Omega, \Sigma, \nu)$ be a $\sigma$-finite measure space, let $a:\R \times \Omega \to \C$ be measurable and bounded, and let $E$ be a spectral measure on $\H$. Then the functions
    \begin{align*}
        \omega &\mapsto \|a(\cdot, \omega)\|_\infty, \qquad        \omega \mapsto \|a(\cdot, \omega)\|_{L^0_\infty(E)}
    \end{align*}
    are measurable.
\end{lem}
\begin{proof}
    Both claims can be proved with the Fubini--Tonelli Theorem or by combining Lemma~\ref{L:PettisIntegral} and Proposition~\ref{P:WeakMeasFunctCalc}. It is vital that $\| \cdot \|_\infty$ is the essential supremum, and that $E$ is a spectral measure on a \textit{separable} Hilbert space as pointed out in~\cite[Remark~4.1.3]{Nikitopoulos2023}.
\end{proof}

\subsection{Abstract Pseudodifferential Calculus}
\label{S:Basics:pseudo}
In this subsection we cover some basics on the pseudodifferential calculus given in Definition~\ref{def:pseudocalc}. See also~\cite{Uuye2011}.

A quick first observation is that the embedding $\H^{s} \hookrightarrow \H^{t}$ is continuous for all $s \geq t$. It follows that $\op^r(\Theta) \subseteq \op^t(\Theta)$ for $r \leq t$.

For $s > 0$, there is a pairing between $\H^s$ and $\H^{-s}$ given by
\[
\langle u, v\rangle_{(\H^s, \H^{-s})} :=\langle  \Theta^s u, \Theta^{-s} v \rangle_\H, \quad u \in \H^s, v\in \H^{-s}.
\]
This pairing identifies $\H^{-s}$ with the (continuous) anti-linear dual space of $\H^s$ and vice-versa.

The space $\H^\infty = \bigcap_{s\in \R} \H^s$ is a Fr\'echet space equipped with the norms $\| \cdot \|_s, s\in \mathbb{R}$. By construction, $\H^\infty \subseteq \H^s$ for any $s \in \R$, and in fact $\H^\infty$ is dense in $\H^s$. Since a subspace of a separable metric space is itself separable it follows that every $\H^s$ admits an orthonormal basis consisting of vectors in $\H^\infty$.

We define $\H^{-\infty}$ as the continuous anti-linear dual space of $\H^\infty$, which can be identified with
\begin{equation}\label{eq:distributions}
    \H^{-\infty} = \bigcup_{s\in \mathbb{R}} \H^{s}.
\end{equation}
This is an $LF$-space, in the sense of~\cite[Chapter~13]{Treves1967}. From this perspective $\H^\infty$ can be interpreted as a Schwartz space and $\H^{-\infty}$ as a space of distributions. 
Given $u \in \H^\infty$ and $v \in \H^{-\infty}$ it follows from~\eqref{eq:distributions} that $v \in \H^{-s}$ for some particular $s \in \R$. It is immediate that $u \in \H^s$, and we have
\[
\langle u, v \rangle_{(\H^{\infty}, \H^{-\infty})} = \langle u, v \rangle_{(\H^{s}, \H^{-s})}.
\]

\begin{prop}\label{P:Interpolation}
    The Sobolev spaces $\H^s$ in Definition~\ref{def:pseudocalc} form an interpolation scale. That is, let $s_0 \leq s_1$, $r_0,r_1\in \R,$ and let $0<\theta<1.$ Set
        \[
            s_\theta := (1-\theta)s_0+\theta s_1,\quad r_\theta = (1-\theta)r_0+\theta r_1.
        \]
        If $T$ is a bounded linear map 
        \[
            T:\H^{s_0}\to \H^{r_0},\quad T\big|_{\H^{s_1}}:\H^{s_1}\to \H^{r_1},
        \]
        then $T|_{\H^{s_\theta}}$ is bounded from $\H^{s_\theta}$ to $\H^{r_{\theta}}$ for every $\theta.$ Moreover we have
        \[
            \|T\|_{\H^{s_\theta}\to \H^{r_{\theta}}} \leq \|T\|_{\H^{s_0}\to \H^{r_0}}^{1-\theta}\|T\|_{\H^{s_1}\to \H^{r_1}}^{\theta}.
        \]
\end{prop}
\begin{proof}
    After identifying $\H^s$ with a weighted $L_2$-space through the spectral theorem, this follows from the Stein--Weiss interpolation theorem for $L_p$-spaces~\cite[Theorem~5.4.1]{BerghLofstrom1976}.
\end{proof}

\subsection{Proof of Theorem~\ref{T:MainMOIConstruction}}
\label{SS:MainProof}

Theorem~\ref{T:MainMOIConstruction} and its proof are heavily inspired by~\cite{ACDS, Peller2006, Peller2016}, see also~\cite{SkripkaTomskova2019}.

\begin{proof}[Proof of Theorem~\ref{T:MainMOIConstruction}]
We have self-adjoint operators $H_0,\ldots, H_n$ on $\H$ with spectral measures $E_j$, and a function $\phi:\R^{n+1}\to\C$ of the form
\[
\phi(\lambda_0,\ldots,\lambda_n)=\int_\Omega a_0(\lambda_0,\omega)\cdots a_n(\lambda_n,\omega) d\nu(\omega),
\]
where $(\Omega,\nu)$ is a finite measure space and the functions $(x, \omega) \mapsto a_j(x,\omega)\langle x \rangle^{-\beta_j}$ are measurable and $E_j \times \nu$-a.e. bounded. Furthermore, $a_j(H_j,\omega)\in\op^{k_j}(\Theta)$, $k_j \in \R$. 

Fix $\omega \in \Omega$ and take $\eta, \xi \in \H^\infty \subseteq \dom a_j(H_j, \omega)$. Then~\cite[Theorem~4.13]{Schmudgen2012} gives that
\[
a_j(H_j, \omega) \xi = \lim_{n\to \infty} a_j(H_j, \omega) \chi_{[-n, n]}(H_j) \xi,
\]
where $\chi_{[-n,n]}$ is the indicator function of the interval $[-n,n]$, because $\operatorname{ess\,sup}_{|\lambda|\leq n} |a_j(\lambda, \omega)| < \infty$. Now Proposition~\ref{P:WeakMeasFunctCalc} gives that
\[
\omega \mapsto \langle \eta, a_j(H_j,\omega) \xi \rangle_{\H^s} =  \lim_{n\to\infty}\langle \Theta^{2s} \eta, a_j(H_j, \omega) \chi_{[-n,n]}(H_j) \xi \rangle_{\H}
\]
is measurable for all $s\in \R$.

Let now $X_i \in \op^{r_i}(\Theta)$, $i \in \{1, \ldots, n\}$. 
 Fix $s \in \R$ and define $s_0, \ldots, s_{2n+1} \in \R$ with
 \[
 s_0:=s, \quad s_{2n+1} := s + \sum_{i=0}^{n} k_i + \sum_{i=1}^n r_i,
 \]
 so that
    the operators $a_j(H_j, \omega)$ and $X_j$ extend to bounded operators
    \begin{align*}
        a_j(H_j, \omega) &\in B(\H^{s_{2j+1}}, \H^{s_{2j}}),\\
        X_j &\in B(\H^{s_{2j}}, \H^{s_{2j-1}}).
    \end{align*}
    By the previous argument, 
    \[
    \omega \mapsto a_j(H_j, \omega) \in B(\H^{s_{2j+1}}, \H^{s_{2j}})
    \]
    is weakly measurable since $\H^{\infty}$ is dense in both $\H^{s_{2j}}$ and $\H^{s_{2j+1}}$.
    
    Using assumption~\eqref{eq:normboundfunctcalc}, i.e.
    \[
    \|a_j(H_j,\omega)\|_{\H \indices{^{s_{2j+1}}}\to\H\indices{^{s_{2j}}}}\leq C_{s,H_j}\|a_j(\cdot,\omega)\|_{L^{\beta_j}_{\infty}(E_j)},
    \]
    we have that
    \begin{align*}
        \int_\Omega \|a_0(H_0, \omega )X_1 a_1(H_1,\omega) &\cdots X_n a_n(H_n, \omega) \|_{\H^{s_{2n+1}} \to \H^{s_0}}d\nu(\omega)\\
        &\lesssim \prod_{j=1}^n \|X_j\|_{\H^{s_{2j}} \to \H^{s_{2j-1}}}  \int_\Omega \prod_{j=0}^n \| a_j(\cdot, \omega)\|_{L^{\beta_j}_\infty(E_j)} d|\nu|(\omega) < \infty,
    \end{align*}
    where Lemma~\ref{L:EssSupMeasurable} ensures the right-hand side is defined. This is a finite quantity since $a_j(x, \omega) \langle x \rangle^{-\beta_j} $ is $E_j \times \nu-a.e.$ bounded and $\nu$ is a finite measure space. Therefore, Lemma~\ref{L:PettisIntegral}
    provides that
    \[
    \int_\Omega a_0(H_0, \omega )X_1 a_1(H_1,\omega) \cdots X_n a_n(H_n, \omega) d\nu(\omega)
    \]
    defines an operator in the weak sense in $B(\H^{s_{2n+1}}, \H^{s_0})$ with
    \begin{align*}
        \bigg \|\int_\Omega a_0(H_0, \omega )X_1 a_1(H_1,\omega) &\cdots X_n a_n(H_n, \omega) d\nu(\omega)\bigg\|_{\H^{s_{2n+1}} \to \H^{s_0}}\\
        &\lesssim \prod_{j=1}^n \|X_j\|_{\H^{s_{2j}} \to \H^{s_{2j-1}}} \int_\Omega \prod_{j=0}^n \| a_j(\cdot, \omega)\|_{L^{\beta_j}_\infty(E_j)} d|\nu|(\omega).
    \end{align*}
    With Pettis' theorem~\cite[Propositions~1.9 and~1.10]{VakhaniaTarieladze1987}, it now follows that for $\psi \in \H^{s_{2n+1}}$, 
    \[
    \omega \mapsto a_0(H_0, \omega )X_1 a_1(H_1,\omega) \cdots X_n a_n(H_n, \omega) \psi \in \H^{s}
    \]
    is Bochner integrable in $\H^s$. This holds in particular for $\psi \in \H^\infty$, and as $s \in \R$ was taken arbitrarily it follows that for $\psi \in \H^\infty$,
    \[
    \int_\Omega  a_0(H_0, \omega )X_1 a_1(H_1,\omega) \cdots X_n a_n(H_n, \omega) \psi d\nu(\omega) \in \H^\infty.
    \]
    It is therefore clear that
    \[
    \int_\Omega  a_0(H_0, \omega )X_1 a_1(H_1,\omega) \cdots X_n a_n(H_n, \omega)  d\nu(\omega) \in \op^{ \sum_j r_j + \sum_j k_j}(\Theta).
    \]
    
That this operator is independent of the chosen representation of $\phi$
\[
\phi(\lambda_0,\ldots,\lambda_n)=\int_\Omega a_0(\lambda_0,\omega)\cdots a_n(\lambda_n,\omega) d\nu(\omega)
\]
follows from the proof of~\cite[Lemma~4.3]{ACDS}. Namely, given $\eta, \xi \in \H^\infty$, it is easy to check that $\theta_{\eta, \xi}:\H^\infty \to \H^\infty$ defined by
\[
\theta_{\eta, \xi} (\psi) := \langle \eta, \psi \rangle_{\H} \xi, \quad \psi \in \H^\infty,
\]
is an element of $\op^{-\infty}(\Theta)$.
The computations in~\cite[Lemma~4.3]{ACDS} give that, for $\eta_k, \xi_k \in \H^\infty$, the integral
\[
\int_\Omega a_0(H_0, \omega) \theta_{\eta_1, \xi_1} a_1(H_1, \omega) \cdots \theta_{\eta_n, \xi_n} a_n(H_n, \omega) d\nu(\omega) \in B(\H)
\]
does not depend on the chosen representation of $\phi$, and so neither does
\[
\int_\Omega a_0(H_0, \omega) \theta_{\eta_1, \xi_1} a_1(H_1, \omega) \cdots \theta_{\eta_n, \xi_n} a_n(H_n, \omega) d\nu(\omega) \bigg|_{\H^\infty} \in \op^{-\infty}(\Theta).
\]
The $SOT$-density of the span of $\{\theta_{\eta, \xi} : \eta, \xi\in \H^\infty\}$ in $B(\H^{s_{2i}}, \H^{s_{2i-1}})$ combined with Lemma~\ref{L:WeaklyMeasComp} concludes the proof.
\end{proof}

\begin{prop}
    The MOI constructed in Theorem~\ref{T:MainMOIConstruction} is linear in its symbol:
    \[
    T_{\alpha\phi + \beta \psi}^{H_0,\ldots, H_n}(X_1, \ldots, X_n) = \alpha T_{\phi}^{H_0,\ldots, H_n}(X_1, \ldots, X_n) + \beta T_{\psi}^{H_0,\ldots, H_n}(X_1, \ldots, X_n), \quad \alpha, \beta \in \C.
    \]
\end{prop}
\begin{proof}
    If both $\phi, \psi: \R^{n+1}\to \C$ have an integral representation of the required form over measure spaces $\Omega$ and $\Sigma$ respectively, then $\alpha \phi + \beta \psi$ can be decomposed appropriately as an integral over the disjoint union $\Omega \sqcup \Sigma$. The assertion then follows by elementary arguments.
\end{proof}

\begin{rem}
    The MOI constructed in Theorem~\ref{T:MainMOIConstruction} is independent of the operator $\Theta$ defining the abstract pseudodifferential calculus in the following sense. If $H_i$ and $X_i$ are operators on $\H$ such that $X_i |_{\H^\infty} \in \op^{r_i}(\Theta)$ and $a_i(H_i, \omega)|_{\H^{\infty}(\Theta)} \in \op^{k_i}(\Theta)$ satisfying the conditions of Theorem~\ref{T:MainMOIConstruction}, then the proof of Theorem~\ref{T:MainMOIConstruction} shows that
    we can define $T_\phi^{H_0,\ldots, H_n}(X_1,\ldots, X_n): \H^{\sum_i r_i + \sum_i k_i} \to \H$ by
    \[
    T_\phi^{H_0, \ldots, H_n}(X_1,\ldots, X_n) \psi = \int_\Omega a_0(H_0, \omega) 
 X_1 a_1(H_1, \omega) \cdots X_n a_n(H_n, \omega) \psi d\nu(\omega) \in \H, \quad \psi \in \H^{\sum_i r_i + \sum_i k_i},
    \]
    which is a map that, apart from the definition of its domain, does not depend on $\Theta$.
\end{rem}

\section{Functional Calculus for Abstract Pseudodifferential Operators}\label{S:FunctCalc}
The construction of multiple operator integrals as abstract pseudodifferential operators now hinges on finding a class of self-adjoint operators $A$ on $\H$ and a class of functions $f$ such that $f(A) \in \op(\Theta)$ with the appropriate estimate of norms~\eqref{eq:normboundfunctcalc}. To accomplish this, and in particular prove Theorem~\ref{T:MainFunctCalc}, we will take $A \in \op^r(\Theta)$, $r \geq 0$, such that $\overline{A}^{r,0}: \H^r \subseteq \H \to \H$ is self-adjoint and study when $f(\overline{A}^{r,0})$ defined via the Borel functional calculus is an operator in $\op(\Theta)$. An obstacle for a naive approach is that $f(\overline{A}^{s+r,s})$ is not easily defined for $s \not =0$, as $\overline{A}^{s+r,s}$ is generally not normal or symmetric. 
As for classical pseudodifferential operators, ellipticity provides the right notion for developing a functional calculus.

\subsection{Elliptic Operators}
\label{SS:EllipticOps}
To prepare the way for a functional calculus on the Sobolev scale, we will show in this subsection that for $A \in \op^r(\Theta)$, $r \geq 0$, $\Theta$-elliptic and symmetric we have that $A$ is self-adjoint with domain $\H^r$. Furthermore, if $A$ is invertible in an appropriate sense, then $A^{-1} \in \op^{-r}(\Theta)$. 

First recall the definition of a $\Theta$-elliptic operator in $\op(\Theta)$.
\begin{defn}
    We say that an operator $A \in \op^r(\Theta)$ is $\Theta$-elliptic if there exists a parametrix for $A$ of order $-r$, that is, there exist operators $P_1, P_2 \in \op^{-r}(\Theta)$ such that
\begin{align*}
    AP &= 1_{\H^{\infty}} + R_1;\\
    PA &= 1_{\H^{\infty}} + R_2,
\end{align*}
where $R_1, R_2 \in \op^{-\infty}(\Theta) = \bigcap_{s \in \R} \op^s(\Theta)$.
\end{defn}
As emphasised in the introduction, the notion of $\Theta$-ellipticity depends on $\Theta$ and on the order $r \in \R$. 

We quickly note that the definition of $\Theta$-ellipticity above is equivalent to the formally weaker condition of $A$ having a right-parametrix $P_1$ and a left-parametrix $P_2$, as it would then follow that $P_1 - P_2 \in \op^{-\infty}(\Theta)$ and hence $P_1$ and $P_2$ are both left- and right-parametrices.

\begin{prop}\label{P:InvertElliptic}
    Let $A \in \op^r(\Theta)$ be $\Theta$-elliptic. If the bounded extension
    \[
    A: \H^{s_0+r} \to \H^{s_0}
    \]
    admits a bounded inverse
    \[
    A^{-1}: \H^{s_0} \to \H^{s_0+r},
    \]
    for a specific $s_0 \in \R$, then $A^{-1}\big|_{\H^\infty} \in \op^{-r}(\Theta)$. Simply writing $A^{-1} = A^{-1}\big|_{\H^\infty}$, we have that
    \[
    A^{-1} A = A A^{-1} = 1_{\H^\infty}.
    \]
\end{prop}
\begin{proof}
    Let $P$ be a parametrix of $A$ and take $x \in \H^\infty$, so that
    \begin{align*}
        A^{-1}x &= (PA - R_2) A^{-1}(AP - R_1) x\\
        &= PAPx - R_2 P x- PR_1 x + R_2 A^{-1} R_1 x.
    \end{align*}
    Observe that for $y \in \H^t$, $t\in \R$, we have $A^{-1}R_1 y \in \H^{s_0+r}$, so that $R_2 A^{-1}R_1 y  \in \H^\infty$.
    Hence,
    \[
    R_2 A^{-1} R_1 \in \op^{-\infty}(\Theta),
    \]
    and therefore
    \begin{align*}
        A^{-1} &= PAP - R_2 P - PR_1 + R_2 A^{-1} R_1 \\
        &\in \op^{-r}(\Theta) + \op^{-\infty}(\Theta) = \op^{-r}(\Theta).
    \end{align*}
\end{proof}

Recall the notion of an asymptotic expansion of pseudodifferential operators~\eqref{eq:pdoexpansion}.

    \begin{lem}[Borel lemma]\label{L:Borel}
        Let $\{A_k\}_{k=0}^\infty$ be a sequence of linear operators from $\H^\infty$ to $\H^\infty$ for which $A_k \in \op^{m_k}(\Theta)$ such that $m_k\downarrow -\infty$
        as $k\to\infty$. There exists a linear operator $A\in \op^{m_0}(\Theta)$ such that
        \begin{equation*}
            A \sim \sum_{k=0}^\infty A_k.
        \end{equation*}
    \end{lem}
    \begin{proof}
        Let $\eta \in C^\infty_c(\R)$ be equal to $1$ in a neighbourhood of zero, and let $\{\varepsilon_k\}_{k=0}^\infty$ be a sequence of positive numbers tending
        to zero in a manner to be determined shortly. Formally we define
        \begin{equation*}
            A := \sum_{k=0}^\infty A_k(1-\eta(\varepsilon_k \Theta)).
        \end{equation*}
        We will prove that $\{\varepsilon_k\}_{k=0}^\infty$ can be chosen such that this series makes sense and $A\in \op^{m_0}(\Theta)$ with the desired asymptotic expansion. 
        
        Let $\xi\in \H^\infty$. Then for every $k \geq 0$ and $n\in \Z$, we have
        \begin{align*}
            \|A_k(1-\eta(\varepsilon_k \Theta))\xi\|_{\H^{n}} &\leq \|A_k\|_{\H^{n+m_k}\to \H^n}\|(1-\eta(\varepsilon_k\Theta))\xi\|_{\H^{n+m_k}}\\
            &\leq \|A_k\|_{\H^{n+m_k}\to \H^n}\|1-\eta(\varepsilon_k\Theta)\|_{\H^{n+m_0}\to \H^{n+m_k}}\|\xi\|_{\H^{n+m_0}}.
        \end{align*}
        Let $a>0$ be a number such that $a < \Theta$. The norm of $1-\eta(\varepsilon_k \Theta)$ from $\H^{n+m_0}$ to $\H^{n+m_k}$ is determined by functional calculus as
        \begin{equation*}
            \sup_{t > a} t^{m_k-m_0} (1-\eta(\varepsilon_k t)) \leq \varepsilon_k^{m_0 - m_k} \sup_{s > 0} s^{m_k - m_0}  (1-\eta(s)) \leq C_\eta\varepsilon_k.
        \end{equation*}
        for some constant $C_\eta$, and for $k$ sufficiently large so that $m_0 - m_k \geq 1$. Now we choose $\varepsilon_k$ sufficiently small such that
        \begin{equation*}
            0 \leq \varepsilon_kC_\eta\max_{|n|\leq k}\{\|A_k\|_{\H^{n+m_k}\to H^{n}}\} < 2^{-k}.
        \end{equation*}
        With this choice of sequence $\{\varepsilon_k\}_{k=0}^\infty$, we have just proved that the series
        \begin{equation*}
            A\xi = \sum_{k=0}^\infty A_k(1-\eta(\varepsilon_k \Theta))\xi
        \end{equation*}
        converges in every $\H^{n}$, and defines a bounded linear operator
        \begin{equation*}
            A:\H^{n+m_0}\to \H^{n},\quad n\in \Z.
        \end{equation*}
        Since this holds for every $n\in \Z$, it follows that $A:\H^\infty\to \H^\infty$ and by interpolation (Proposition~\ref{P:Interpolation}) $A \in \op^{m_0}(\Theta)$.
        Note that with this fixed choice of $\{\varepsilon_k\}_{k=0}^\infty$, we have proved the stronger result that the ``tail" of $A$
        \begin{equation*}
            \sum_{k=N+1}^\infty A_k(1-\eta(\varepsilon_k \Theta))
        \end{equation*}
        converges in every $\H^s$ and defines a linear operator in $\op^{m_{N+1}}(\Theta)$.
        
        Now we prove that $A$ has the desired asymptotic expansion.
        For every $N>0$ we have
        \begin{equation*}
            A-\sum_{k=0}^N A_k = -\sum_{k=0}^N A_k\eta(\varepsilon_k \Theta) + \sum_{k=N+1}^\infty A_k(1-\eta(\varepsilon_k \Theta))
        \end{equation*} 
        Since $\eta$ is compactly supported, it is easy to see that the first summand has order $-\infty$ for every $N\geq 0$, and the second summand
        has order at most $m_{N+1}$ due to the result just proved. 
    \end{proof}

\begin{cor}\label{C:EllipticParametrix}
    Suppose that $A \in \op^r(\Theta)$ has an inverse $B\in \op^{-r}(\Theta)$ up to order $-1$. That is, 
        \begin{align*}
            AB = 1_{\H^\infty}+R_1;\\
            BA = 1_{\H^\infty}+R_2
        \end{align*} 
        where $R_1, R_2$ have order $-1$. Then $A$ is $\Theta$-elliptic.
\end{cor}
\begin{proof}
    Since $R_1^j$ has order at most $-j$, we can use the Borel lemma to construct an operator $B'$ such that
    \begin{equation*}
            B' \sim \sum_{k=0}^\infty BR_1^j.
        \end{equation*}
        Then $AB' -1$ has order $-\infty$. Similarly we can construct a left inverse.
\end{proof}

\begin{prop}\label{P:EllipticClosed}
    Let $A \in \op^r(\Theta)$ be $\Theta$-elliptic of order $r \geq 0$. Then the unbounded operator
    \[
    A: \H^{s+r}\subseteq \H^s \to \H^s
    \]
    is closed for each $s \in \R$.
\end{prop}
\begin{proof}
    Define the graph norm of $A$ on $\dom(A) = \H^{s+r}$ as 
    \[
    \| x \|_{G(A)} := \|Ax \|_{s} + \| x\|_s, \quad x \in \H^{s+r}.
    \]
    By definition, $A$ is a closed operator if and only if $\dom(A)$ is complete with respect to this graph norm. We will show that for $\Theta$-elliptic operators, the graph norm is equivalent to $\| \cdot \|_{s+r}$, which immediately implies that $\H^{s+r}$ is complete with respect to the graph norm. First, we have that
    \begin{align*}
        \|Ax \|_s + \|x\|_s & \leq \|A\|_{\H^{s+r}\to \H^s} \|x\|_{s+r} + \|\Theta^{-r}\|_{\H^s \to \H^s} \|\Theta^r x \|_s\\
        & \lesssim \| x\|_{s+r}.
    \end{align*}
    Next, let $P$ be a parametrix for $A$, and take $x \in \H^\infty$ so that
    \begin{align*}
        \| x \|_{s+r} & \leq \| PA x \|_{s+r} + \| R_2 x\|_{s+r}\\
        & \leq \| P \|_{\H^s \to \H^{s+r}} \| Ax \|_{s} + \| R_2 \|_{\H^s \to \H^{s+r}} \| x \|_s\\
        &\lesssim \|Ax \|_s + \|x\|_s.
    \end{align*}
    The assertion of the proposition is now immediate.
\end{proof}

The $\Theta$-elliptic operators have a property which is often called elliptic regularity, or maximal subellipticity.
\begin{prop}\label{P:ellipticreg}
    Let $A \in \op^r(\Theta)$ be a $\Theta$-elliptic operator. If $x \in \H^{-\infty}$ is such that $Ax \in \H^s$ for an $s \in \R$, then $x \in \H^{s+r}$.
\end{prop}
\begin{proof}
    Take $x \in \H^{-\infty}$, and suppose that $Ax \in \H^{s}$. Then $ PAx\in \H^{s+r}$,
    which implies that
    \[
    x = PAx - R_2x \in \H^{s+r}.
    \]
\end{proof}

Finally we will now show that if $A \in \op^r(\Theta)$, $r \geq 0$ is $\Theta$-elliptic and symmetric, then $A$ is self-adjoint with domain $\H^r$. 
\begin{defn}
    For $A \in \op^r(\Theta)$, we define the adjoint $A^\dag:\H^{-\infty} \to \H^{-\infty}$ by
    \[
    \langle Au, v \rangle_{(\H^\infty, \H^{-\infty})} = \langle u, A^\dag v \rangle_{(\H^\infty, \H^{-\infty})} \quad u \in \H^\infty, v\in \H^{-\infty}.
    \]
\end{defn}
It is shown in Appendix~\ref{S:AppEllipticAdj} that $A^\dag|_{\H^\infty} \in \op^r(\Theta)$, and for any $s\in \R$, we have that
\begin{equation}\label{eq:DaggerAdj}
    \langle A u, v \rangle_{(\H^s, \H^{-s})} = \langle u, A^\dag v \rangle_{(\H^{s+r}, \H^{-s-r})}, \quad u \in \H^{s+r}, v\in \H^{-s}.
\end{equation}
Furthermore, for $r \geq 0$ Proposition~\ref{P:symmetric} gives that $A = A^\dag$ if and only if $\overline{A}^{r,0}: \H^r \subseteq \H \to \H$ is symmetric.

\begin{prop}\label{P:EllipticSelfAdj}
    Let $A \in \op^r(\Theta)$, $r \geq 0$, be a $\Theta$-elliptic and symmetric operator. Then $A$
    is self-adjoint with domain $\H^r$.
\end{prop}
\begin{proof}
    We denote the Hermitian adjoint of the closed operator $\overline{A}^{r,0}: \H^{r} \subseteq \H^0 \to \H^0$ by $A^{*_0}:= \big(\overline{A}^{r,0} \big)^*$,
    \[
    A^{*_0}: \dom(A^{*_0}) \subseteq \H^0 \to \H^0.
    \]
    In order to prove that $A$ is self-adjoint, we need to show that $\dom(A^{*_0}) = \H^r$.
    Recall that, by definition,
    \begin{align*}
        \dom(A^{*_0}) := \{u \in \H^0 : \exists v \in \H^0 \text{ such that } \forall \phi \in \H^r \ \langle u, A\phi \rangle_{\H^0}=\langle v, \phi \rangle_{\H^0} \}.
    \end{align*}
    If $u,v \in \H^0$ and $\phi \in \H^r$, then by~\eqref{eq:DaggerAdj} and Proposition~\ref{P:symmetric},
    \begin{align*}
        \langle u, A\phi \rangle_{\H^0} &=\langle A^\dag u, \phi \rangle_{(\H^{-r}, \H^{r})}\\
        &= \langle A u, \phi \rangle_{(\H^{-r}, \H^{r})};\\
        \langle v, \phi \rangle_{\H^0} &= \langle v, \phi \rangle_{(\H^{-r}, \H^{r})}.
    \end{align*}
    Since $\H^{r}$ separates the points of $\H^{-r}$, we have that for $u, v \in \H^0$,
    \[
    \langle A u, \phi \rangle_{(\H^{-r}, \H^{r})} = \langle v, \phi \rangle_{(\H^{-r}, \H^{r})}, \quad \forall \phi \in \H^{r},
    \]
    if and only if
    \[
    Au = v \in \H^{-r}.
    \]
    Hence, 
    \begin{align*}
        \dom(A^{*_0}) &= \{u \in \H^0 : \exists v \in \H^0 \text{ such that }  A u = v \} \\
        &= \{u \in \H^0 :  A u \in \H^0\}.
    \end{align*}
    By elliptic regularity (Proposition~\ref{P:ellipticreg}) it follows that $\dom(A^{*_0}) = \H^r,$
    completing the proof.
\end{proof}

\subsection{Functional Calculus for Elliptic Operators}
\label{SS:FunctCalc}
\begin{prop}\label{P:ReplaceTheta}
    Let $A\in \op^r(\Theta)$, $r>0$ be $\Theta$-elliptic and symmetric. Then $\langle A \rangle^{\frac{1}{r}}$ extends to an invertible positive self-adjoint operator on $\H$, and 
    \[
    \dom \Theta^s = \dom \big(\langle A \rangle^{\frac{1}{r}}\big)^s.
    \]
    The norms $\| \Theta^s \xi \|_{\H}$ and $\|\langle A \rangle^{\frac{s}{r}} \xi \|_{\H}$ are equivalent on this subspace of $\H$. Therefore, $\Theta$ and $\langle A \rangle^{\frac{1}{r}}$ define the same Sobolev scale
    \[
    \H^s(\Theta) = \H^s\big(\langle A \rangle^{\frac{1}{r}}\big),
    \]
    and we have
    \[
    \op^t(\Theta)  = \op^t \big(\langle A \rangle^{\frac{1}{r}}\big).
    \]
\end{prop}
\begin{proof}
The first statement follows from Proposition~\ref{P:EllipticSelfAdj}. For the remaining statements, it suffices to prove that 
    \[
    (1+A^2)^{\alpha} \Theta^{-2\alpha r}
    \]
    extends to a bounded operator on $\H$ for all $\alpha \in \R$, as this would imply for $\xi \in \H^\infty$,
    \[
    \|\langle A \rangle^{\frac{s}{r}} \xi \|_{\H} \leq \| (1+A^2)^{\frac{s}{2r}} \Theta^{-s}\|_\infty \| \Theta^s \xi\|_{\H} \lesssim \| \Theta^s \xi\|_{\H},
    \]
    and an analogous estimate in the other direction.

Let $P$ be a parametrix for $A$, so that $AP = 1 + R$ with $R \in \op^{-\infty}(\Theta)$. Since
    \[
    (1+A^2)P^2 = P^2 + A(1+R)P = 1 + P^2 + R + ARP
    \]
    and $P^2 + R + ARP \in \op^{-2r}(\Theta)$ and similarly for $P^2(1+A^2)$, it follows that the operator $1+A^2$ is also $\Theta$-elliptic due to Corollary~\ref{C:EllipticParametrix}. Since $A$ is self-adjoint with domain $\H^r$, applying Proposition~\ref{P:InvertElliptic} gives that  $(1+A^2)^{-1} \in \op^{-2r}(\Theta)$. We therefore have $(1+A^2)^k \in \op^{2kr}(\Theta)$, $k \in \Z$. This in turn gives that $\H^\infty \subseteq \dom (1+A^2)^z$ for any $z\in \C$.
    
    We use the Hadamard three-line theorem, so define the function
    \[
    F(z) := \langle x, (1+A^2)^{mz} \Theta^{-2 mz r} y \rangle_{\H}, \quad z \in \C,
    \]
    where $m \in \Z$ and $x, y \in \H^\infty$ are fixed.     
    Let $\{e_n\}_{n\in \N}\subseteq \H^\infty$ be an orthonormal basis of $\H$, then
    \begin{equation}\label{E:basis_expansion_of_F}
    F(z) = \sum_{n=1}^\infty \langle x, (1+A^2)^{mz} e_n \rangle_\H \langle e_n , \Theta^{-2mzr} y\rangle_{\H}.
    \end{equation}
    Using the dominated convergence theorem, it can be seen that $z \mapsto  \langle x, (1+A^2)^{mz} e_n \rangle_\H$ and $z \mapsto \langle e_n , \Theta^{-2mzr} y\rangle_{\H}$ are continuous maps. Applying the Cauchy--Schwarz inequality to the series \eqref{E:basis_expansion_of_F} yields
    \[
        \sum_{n=1}^\infty |\langle x,(1+A^2)^{mz}e_n\rangle_{\H}||\langle e_n,\Theta^{-2mzr}y\rangle_{\H}| \leq \|(1+A^2)^{m\overline{z}}x\|_{\H}\|\Theta^{-2mzr}y\|_{\H},
    \]
    which is uniformly bounded on compact subsets of $\C$ due to the continuity of the right-hand side. 
    We can therefore apply the dominated convergence theorem again to deduce that $F(z)$ is a continuous function itself. 
    Furthermore, this uniform boundedness yields through Fubini's theorem that if $\gamma$ is a closed loop
    in $\mathbb{C}$ then
    \[
        \int_{\gamma} F(z)\,dz = \sum_{n=1}^\infty \int_{\gamma}\langle x,(1+A^2)^{mz}e_n\rangle_{\H}\langle e_n,\Theta^{-2mzr}y\rangle_{\H}\,dz.
    \]
    Using Fubini's theorem once more,
    we have that \[
    \int_\gamma F(z) dz = \sum_{n=1}^\infty  \int_{\sigma(1+A^2)}\int_{\sigma(\Theta^{-2r})} \int_{\gamma} (\lambda \mu)^{mz} dz \langle x,dE^{1+A^2}e_n\rangle_{\H} \langle e_n,dE^{\Theta^{-2r}}y\rangle_{\H} = 0,
    \] so that we can conclude by Morera's theorem that $F(z)$ is holomorphic.    
    
    Since $1+A^2$ and $\Theta$ are positive operators and
    \[
    \sup_{x > 0} |x^{it}| = 1,
    \]
    it follows from the Borel functional calculus that for $s \in \R,$
    \begin{align*}
        |F(is)| &= |\langle  (1+A^2)^{-ims} x, \Theta^{-2 imsr} y \rangle_{\H}|\\
        &\leq \| x \|_{\H} \|y\|_{\H}.
    \end{align*} 
    Likewise,
    \begin{align*}
        |F(1+is)| &= |\langle  (1+A^2)^{-ims} x, (1+A^2)^m \Theta^{-2mr} \Theta^{-2imsr} y \rangle_{\H}|\\
        &\leq \| (1+A^2)^m \Theta^{-2mr} \|_{\H \to \H} \|x \|_{\H} \|y\|_{\H},
    \end{align*}
    which we know to be finite since $m$ is an integer.

    The Hadamard three-line theorem (see e.g. \cite[Lemma 1.1.2]{BerghLofstrom1976}) now gives that for $\alpha \in (0,1)$,
    \begin{align*}
        |F(\alpha)| &\leq \max_{s\in \R} |F(\alpha + is)|\\
        &\leq \big( \max_{s\in \R} |F(is)|\big)^{(1-\alpha)} \big( \max_{s\in \R} |F(1 + is)|\big)^\alpha\\
        &\leq  \| (1+A^2)^m \Theta^{-2mr} \|_{\H \to \H}^\alpha \|x\|_{\H} \|y\|_{\H}.
    \end{align*}
    Hence, with $\alpha \in (0,1)$ and $m \in \Z$,
    \begin{align*}
        \|(1+A^2)^{m\alpha} \Theta^{-2m \alpha r}\|_{\H\to \H} \leq  \| (1+A^2)^m \Theta^{-2mr} \|_{\H \to \H}^\alpha,
    \end{align*}
    which proves that $(1+A^2)^{m\alpha} \Theta^{-2m\alpha r}$ extends to a bounded operator on $\H$ for all $m\in \Z$, $\alpha \in [0,1]$.
\end{proof}

\begin{proof}[Proof of Theorem~\ref{T:MainFunctCalc}]
    Let $A \in \op^r(\Theta)$, $r>0$ be $\Theta$-elliptic and symmetric. Using Proposition~\ref{P:ReplaceTheta}, we replace $\Theta$ by $\langle A \rangle^{\frac{1}{r}}$ so that $A \in \op^r(\langle A \rangle^{\frac{1}{r}})$ is $\Theta$-elliptic and symmetric. By Proposition~\ref{P:EllipticSelfAdj}, the operator
    \[
    \overline{A}^{r,0}: \H^r(\langle A \rangle^{\frac{1}{r}}) \subseteq \H \to \H
    \]
    is self-adjoint; we denote its spectral measure by $E$. Then for $f \in L^\beta_\infty(E)$, using Borel functional calculus to define $f(\overline{A}^{r,0})$, we have
    \begin{align*}
        \|\langle A \rangle^{\frac{s}{r}} f(\overline{A}^{r,0})\langle A \rangle^{-\frac{s}{r} - \beta}\|_{\H \to \H} = \| f(\overline{A}^{r,0})\langle A \rangle^{- \beta}\|_{\H \to \H} = \|f\|_{L^\beta_\infty(E)}<\infty,
    \end{align*}
    which shows that $f(A) := f(\overline{A}^{r,0})|_{\H^\infty} \in \op^{r\beta}(\langle A \rangle^{\frac{1}{r}})$. Converting this estimate back into an estimate on the spaces $\H^s(\Theta)$ introduces the constant $C_{s,A}$.
\end{proof}

Theorem~\ref{T:MainFunctCalc} has a converse in the following sense. If $A \in \op^r(\Theta)$, $r>0$ is an arbitrary $\Theta$-elliptic symmetric operator and if $f: \R \to \C$ is such that $f(A) \in \op^{\beta r}(\Theta)$, then the proof of Proposition~\ref{P:ReplaceTheta} gives that $f(A) (1+A^2)^{-\beta/2}$ is a bounded operator on $\H$. This happens if and only if $f \langle x \rangle^{-\beta} \in L^0_\infty(E)$~\cite[Theorem~5.9]{Schmudgen2012}, i.e. $f\in L^\beta_\infty(E)$.

\begin{cor}\label{C:IndepNorm}
    If $A\in \op^r(\Theta), r>0$ is symmetric and $\Theta$-elliptic and if $f :\R \to \C$ is a bounded Borel measurable function, then for any $t\in \R$ we have
    \[
    \|f(tA)\|_{\H^s \to \H^s} \leq C_{s,A} \|f\|_{\infty},
    \]
    independent of $t\in \R$.
\end{cor}
\begin{proof}
    Like in the proof Theorem~\ref{T:MainFunctCalc}, we have
    \begin{align*}
        \|f(tA)\|_{\H^s(\langle A \rangle^{\frac{1}{r}}) \to \H^s(\langle A \rangle^{\frac{1}{r}})} & = \|f_t\|_{L^0_\infty(E)} \leq \|f\|_\infty,
    \end{align*}
    where we wrote $f_t(x) := f(tx)$.
\end{proof}

The functional calculus constructed in Theorem~\ref{T:MainFunctCalc} can easily be extended to a larger class of operators. For example, on $\R^d$ with $\Theta = (1-\Delta)^{1/2}$ we have that $i\frac{d}{dx}$ is not $ (1-\Delta)^{1/2}$-elliptic, but it does commute strongly with a symmetric elliptic operator. The following proposition shows that a functional calculus for $i\frac{d}{dx}$ does exist in $\op(\Theta)$ for this reason.
\begin{prop}\label{P:BigFunctCalc}
    Let $A$ be a self-adjoint operator on $\H$ with spectral measure $E$. If there exists a $\Theta$-elliptic symmetric operator $H \in \op^h(\Theta)$, $h>0$ such that $A$ strongly commutes with $\overline{H}^{h,0}: \H^h \subseteq \H \to \H$, then for $f \in L^0_\infty(E)$ we have that $f(A) \in \op^0(\Theta)$. If $A \in \op^r(\Theta)$ itself for some $r \in \R$, we have that $f(A) \in \op^{\beta r}(\Theta)$ for $f \in L^\beta_\infty(E)$, $\beta \geq 0$.
\end{prop}
\begin{proof}
    In light of Proposition~\ref{P:ReplaceTheta}, we can assume without loss of generality that $H = \Theta$. If $f \in L^0_\infty(E)$, then $f(A):\H \to \H$ is a bounded operator, and for $\xi \in \H^\infty$ we have
    \[
    \| f(A) \xi\|_{\H^k} = \|\Theta^k f(A) \xi\|_{\H} \leq \|f(A)\|_\infty \|\xi\|_{\H^k}, \quad k \in \Z, 
    \]
    which shows that $f(A) \in \op^0(\Theta)$ through interpolation (Proposition~\ref{P:Interpolation}). The second part of the proposition is proved similarly, after the observation that the Hadamard three-line argument in the proof of Proposition~\ref{P:ReplaceTheta} goes through for $A$ if $m \in \Z_{\geq 0}$, i.e. $(1+A^2)^\alpha \Theta^{-2\alpha r}$ is bounded for $\alpha \geq 0$. 
\end{proof}

\section{Divided differences}\label{S:FunctSpac}
The condition that appears in Theorem~\ref{T:MainMOIConstruction} on the symbol $\phi: \R^{n+1}\to \C$ needs to be analysed more closely in order to prepare this multiple operator integral construction for practical applications. The main result of this section is Lemma~\ref{L:DivDifT}, which gives that for $f \in T^\beta(\R)$ the divided difference $f^{[n]}$ has an integral representation satisfying the conditions of Theorem~\ref{T:MainMOIConstruction}.

First of all, for functions $\phi:\R^{n+1} \to \C$ it is an equivalent condition to admit a representation
\[
\phi(\lambda_0,\ldots,\lambda_n)=\int_\Omega a_0(\lambda_0,\omega)\cdots a_n(\lambda_n,\omega) d\nu(\omega)
\]
for a finite complex measure space (technically: of finite variation) $(\Omega,\nu)$ and measurable functions $a_j:\R\times\Omega\to\C$ such $(x,\omega)\mapsto a_j(x,\omega)\langle x \rangle^{-\beta_j}$ is $E_j \times \nu$-a.e. bounded for $\beta_j \in \R$,
or a representation
\[
\phi(\lambda_0,\ldots,\lambda_n)=\int_\Sigma b_0(\lambda_0,\sigma)\cdots b_n(\lambda_n,\sigma) d\mu(\sigma),
\]
where $\Sigma$ is a $\sigma$-finite measure space, $b_j: \R \times \Omega \to \C$ measurable, and
\[
\int_\Sigma \|b_0(\cdot, \sigma)\|_{L^{\beta_0}_\infty(E_0)} \cdots \|b_n(\cdot, \sigma)\|_{L^{\beta_n}_\infty(E_n)}d|\mu|(\sigma)<\infty.
\]
Namely, given the second representation, a representation of the first type can be obtained by putting~\cite[p.48]{SkripkaTomskova2019}
\[
a_j(\lambda_j, \sigma) := \frac{b_j(\lambda_j, \sigma)}{\|b_j(\cdot, \sigma)\|_{L_\infty^{\beta_j}(E_j)}}, \quad \nu(\sigma):= \|b_0(\cdot, \sigma)\|_{L^{\beta_0}_\infty(E_0)} \cdots \|b_n(\cdot, \sigma)\|_{L^{\beta_n}_\infty(E_n)}\mu(\sigma).
\]

Mainly for notational purposes, we introduce the following definition. It is inspired by the integral projective tensor product appearing in~\cite{Peller2006}, a precise study can be found in~\cite{Nikitopoulos2023}.

\begin{defn}\label{def:ProjIntBox}
    Let $\Gamma_0, \ldots, \Gamma_n$ be function spaces of bounded measurable functions $\R \to \R$ equipped with (semi)norms $\|\cdot \|_{\Gamma_i, k}$.
    We define $\Gamma_0 \boxtimes_i \cdots \boxtimes_i \Gamma_n$ as the set of functions $\phi: \mathbb{R}^{n+1} \to \mathbb{C}$ for which there exists a decomposition
    \begin{equation}\label{eq:boxrep}
        \phi(\lambda_0, \ldots, \lambda_n) = \int_\Omega a_0(\lambda_0, \omega) \cdots a_n(\lambda_n, \omega) d\nu(\omega)
    \end{equation}
    where $(\Omega, \nu)$ is a $\sigma$-finite measure space, $a_i : \mathbb{R} \times \Omega \to \mathbb{C}$ is measurable, $a_i(\cdot, \omega) \in \Gamma_i$, the functions $\omega \mapsto \|a_i(\cdot, \omega)\|_{\Gamma_i, k}$ are measurable for each $i$ and $k$, and
    \begin{equation}\label{eq:boxnorm}
        \int_\Omega \|a_0(\cdot, \omega)\|_{\Gamma_0, k_0} \cdots \|a_n(\cdot, \omega)\|_{\Gamma_n, k_n} d|\nu|(\omega) < \infty.
    \end{equation}
    We define the seminorm
    \[
    \| \phi \|_{\Gamma_0 \boxtimes_i \cdots \boxtimes_i \Gamma_n, k_0, \ldots, k_n}
    \]
    to be the infimum of the quantity~\eqref{eq:boxnorm} over all representations~\eqref{eq:boxrep}.
\end{defn}

\begin{rem}
We have that 
\[
\Gamma_0 \otimes \cdots\otimes \Gamma_n \subseteq \Gamma_0 \boxtimes_i \cdots \boxtimes_i \Gamma_n,
\]
where $\otimes$ denotes the algebraic tensor product for topological vector spaces.
\end{rem}

Note that due to Lemma~\ref{L:EssSupMeasurable}, if $a: \R \times \Omega \to \C$ is measurable and $a(\cdot, \omega) \in S^\beta(\R)$, we have that $\omega \mapsto \|a(\cdot, \omega)\|_{S^\beta(\R), k}$ is measurable, and the same claim holds for $T^\beta(\R)$ and $L_\infty^\beta(E)$. Hence the construction in Definition~\ref{def:ProjIntBox} can be applied without this extra assumption.

For $L_\infty^\beta(E)$ with $E$ a spectral measure, this gives the integral projective tensor product
\[
L_\infty^{\beta_0}(E_0) \boxtimes_i \cdots \boxtimes_i L_\infty^{\beta_n}(E_n) = L_\infty^{\beta_0}(E_0) \hat{\otimes}_i \cdots \hat{\otimes}_i L_\infty^{\beta_n}(E_n),
\]
which appears in particular for $\beta_0 = \cdots = \beta_n = 0$ in the works by Peller~\cite[pp.6, 7]{Peller2006}.

We refrain from answering the question whether the space in Definition~\ref{def:ProjIntBox} is in general the completion of the algebraic tensor product under the given seminorms. 

Recall the previously defined seminorms
\begin{align*}
\|f\|_{S^\beta(\R), k} &= \sup_{x\in \R} |f^{(k)}(x)| \langle x \rangle^{k-\beta} < \infty, \quad k \in \Z_{\geq 0},\\
    \|f\|_{T^\beta(\R), k}&= \int_{\R} |f^{(k)}(x)|\langle x\rangle^{k-\beta-1}dx < \infty, \quad k \in \Z_{\geq 0},
    \end{align*}
from Definition~\ref{def:SbetaTbeta}

\begin{rem}\label{rem:SipvL}
    Observe that
    \[
    S^{\beta_0}(\R) \boxtimes_i \cdots \boxtimes_i S^{\beta_n}(\R) \subseteq  L_\infty^{\beta_0}(E_0) \hat{\otimes}_i \cdots \hat{\otimes}_i L_\infty^{\beta_n}(E_n)
    \]
    no matter what spectral measures $E_i$ are taken. 
    By Theorem~\ref{T:MainMOIConstruction} we have a well-defined multiple operator integral precisely for symbols in the latter space.
\end{rem}

\begin{prop}\label{P:BoxProdProd}
If $\phi \in S^{\alpha_0}(\R) \boxtimes_i \cdots \boxtimes_i S^{\alpha_n}(\R)$ and $\psi \in S^{\beta_0}(\R) \boxtimes_i \cdots \boxtimes_i S^{\beta_n}(\R)$, then
\[
\Phi(\lambda_0, \ldots, \lambda_n):= \phi(\lambda_0, \ldots, \lambda_n) \psi(\lambda_0, \ldots, \lambda_n)\in S^{\alpha_0 + \beta_0}(\R) \boxtimes_i \cdots \boxtimes_i S^{\alpha_n + \beta_n}(\R).
\]
An analogous statement holds for the spaces $L^\beta_\infty(E)$.
\end{prop}
\begin{proof}
    According to Definition~\ref{def:ProjIntBox}, we can find $\sigma$-finite measure spaces $(\Omega, \nu)$, $(\Sigma, \mu)$ and measurable functions $a_i: \R \times \Omega \to \C$, $b_i: \R \times \Sigma \to \C$ such that
    \begin{align*}
        \phi(\lambda_0, \ldots, \lambda_n) &= \int_{\Omega} a_0(\lambda_0, \omega) \cdots a_n(\lambda_n, \omega) d\nu(\omega);\\
        \psi(\lambda_0, \ldots, \lambda_n) &= \int_{\Sigma} b_0(\lambda_0 ,\sigma) \cdots b_n(\lambda_n, \sigma) d\mu(\sigma).
    \end{align*}
As observed above, the maps
    \[
    \omega \mapsto \| a_i(\cdot, \omega)\|_{S^{\alpha_i}(\R), k}
    \]
    are measurable, similarly for the functions $b_i$.
    
    Using Tonelli's theorem,
    \begin{align*}
        \int_{\Omega\times \Sigma} \big|a_0(\lambda_0, \omega)b_0(\lambda_0, \sigma)&\cdots a_n(\lambda_n, \omega)b_n(\lambda_n, \sigma)\big| d(\nu\times \mu)(\omega, \sigma)\\
        &\leq \langle \lambda_0 \rangle^{\alpha_0+\beta_0} \cdots \langle \lambda_n \rangle^{\alpha_n + \beta_n} \int_{\Omega}  \| a_0(\cdot, \omega)\|_{S^{\alpha_0},0} \cdots \|a_n(\cdot, \omega)\|_{S^{\alpha_n}, 0} d\nu(\omega) \\
        &\quad \times \int_{\Sigma}\|b_0(\cdot, \sigma)\|_{S^{\beta_0}, 0} \cdots \|b_n(\cdot, \sigma)\|_{S^{\beta_n},0} d\mu( \sigma)<\infty.
    \end{align*}
    Hence, by Fubini's theorem
    \[
    \Phi(\lambda_0, \ldots, \lambda_n) = \int_{\Omega \times \Sigma} a_0(\lambda_0, \omega)b_0(\lambda_0, \sigma) \cdots a_n(\lambda_n, \omega) b_n(\lambda_n, \sigma) d(\nu \times \mu)(\omega, \sigma).
    \]
    The fact that $\Phi \in S^{\alpha_0+\beta_0}(\R) \boxtimes_i \cdots \boxtimes_i S^{\alpha_n + \beta_n}(\R)$ now follows from the computation
    \[
    \| a_k(\cdot, \omega) b_k(\cdot, \sigma)\|_{S^{\alpha_k + \beta_k}, m} \leq \sum_{j=0}^m \binom{m}{j} \|a_k(\cdot, \omega)\|_{S^{\alpha_k}, j} \| b_k(\cdot, \sigma)\|_{S^{\beta_k}, m-j}.
    \]
\end{proof}

We now use the construction of an almost analytic extension to provide an explicit integral representation for $f\in T^\beta(\R)$. The technique was introduced by H\"ormander~\cite{Hormander1969, Hormander1970}, and subsequently used in various contexts by many authors. For detailed notes on the historical origins, see~\cite{Hormander1969}. In Appendix~\ref{S:AppFunctCalc0} we use it to develop a functional calculus for $\op^0(\Theta)$, based on the results of~\cite{Davies1995a}. 

For the details of this construction we follow Davies~\cite{Davies1995a}\cite[Section~2.2]{Davies1995c}.
\begin{defn}[\cite{Davies1995a}]\label{def:AlmostAnalyticExtension1}
    Let $f \in C^N(\mathbb{R})$, $N \in \N$. We define an extension $\tilde{f}: \mathbb{C} \to \mathbb{C}$ by
    \[
    \tilde{f}(x+iy) := \tau(y/\langle x \rangle ) \sum_{r=0}^N f^{(r)}(x) \frac{(iy)^r}{r!},
    \]
    where $\tau:\mathbb{R} \to \mathbb{R}$ is a smooth bump function with $\tau(s) = 0$ for $|s|>2$, $\tau(s) =1$ for $|s|<1$. If $f \in C_c^\infty(\R)$, we have
    \[
    f(x) = -\frac{1}{\pi} \int_\mathbb{C} \left( \frac{\partial \tilde{f}}{\partial \overline{z}}(z)\right) (z-x)^{-1} dz, \quad x \in \mathbb{R},
    \]
    independent of the choice of $\tau$ and $N$. We refer to $\tilde{f}$ as an \textit{almost analytic extension} of $f$.
\end{defn}

\begin{lem}\label{L:DivDifT}
\begin{enumerate}
    \item Take $n \in \Z_{\geq 0}$, let $\alpha$ be some real number with $-1 \leq \alpha \leq n$, and consider any collection of real numbers $-1 \leq \beta_0, \ldots, \beta_n \leq 0$ such that $\sum \beta_j = \alpha - n$. For each $k_0, \ldots, k_n \geq 0$, and $f \in C^{n+\sum_{j=0}^n k_j+2}(\R)$, we have
\begin{align*}
\|f^{[n]}(\lambda_0, \ldots, \lambda_n)\|_{S^{\beta_0}(\R) \boxtimes_i \cdots \boxtimes_i S^{\beta_n}(\R), k_0, \ldots, k_n} \lesssim \sum_{r=0}^{n+\sum_{j=0}^n k_j+2} \| f \|_{T^\alpha(\R), r}.
\end{align*}
In particular,
    \[
f \in T^\alpha(\R) \Rightarrow f^{[n]} \in S^{\beta_0}(\R) \boxtimes_i \cdots \boxtimes_i S^{\beta_n}(\R),
\]

\item Let $\alpha \leq -1$, and $f\in C^{n+\sum_{j=0}^n k_j+2}(\R)$. Then $f^{[n]}$ can be decomposed into a finite sum of functions $\phi$ for which
\begin{align*}
\|\phi\|_{S^{\beta_0}(\R) \boxtimes_i \cdots \boxtimes_i S^{\beta_n}(\R), k_0, \ldots, k_n} \lesssim \sum_{r=0}^{n+\sum_{j=0}^n k_j+2} \| f \|_{T^\alpha(\R), r},
\end{align*}
with $\beta_0,\ldots, \beta_n \leq 0$ and $\sum_j \beta_j = \alpha - n$. In particular,
\[
f \in T^\alpha(\R) \Rightarrow f^{[n]} \in \sum_{\substack{\beta_0, \ldots, \beta_n \leq 0 \\ \sum \beta_j= \alpha -n}}S^{\beta_0}(\R) \boxtimes_i \cdots \boxtimes_i S^{\beta_n}(\R).
\]

\item Let $\alpha \geq n$, and $f\in C^{n+\sum_{j=0}^n k_j+2}(\R)$. Then $f^{[n]}$ can be decomposed into a finite sum of functions $\phi$ for which
\begin{align*}
\|\phi\|_{S^{\beta_0}(\R) \boxtimes_i \cdots \boxtimes_i S^{\beta_n}(\R), k_0, \ldots, k_n} \lesssim \sum_{r=0}^{n+\sum_{j=0}^n k_j+2} \| f \|_{T^\alpha(\R), r},
\end{align*}
where $\sum_j \beta_j = \alpha -n $. 
In particular,
\[
f \in T^\alpha(\R) \Rightarrow f^{[n]} \in \sum_{ \sum \beta_j= \alpha -n}S^{\beta_0}(\R) \boxtimes_i \cdots \boxtimes_i S^{\beta_n}(\R).
\]

\end{enumerate}
\end{lem}
\begin{proof}
\begin{enumerate}
    \item 
 Take $f\in C^{N+1}(\R)$ 
with almost analytic extension $\tilde{f}$. Directly from Definition~\ref{def:AlmostAnalyticExtension1}, writing $\sigma(z) := \tau(\frac{\Im(z)}{\langle\Re(z)\rangle})$, it follows that (cf. \cite[Section~2.2]{Davies1995b})
\begin{align*}
    \frac{\partial \tilde{f}}{\partial \overline{z}} = \frac{1}{2} \bigg(\sum_{r=0}^N f^{(r)}(\Re(z)) \frac{(i\Im(z))^r}{r!} \bigg)(\sigma_x(z) + i\sigma_y(z) ) + \frac{1}{2} f^{(N+1)}(\Re(z)) \frac{(i\Im(z))^N}{n!}\sigma(z).
\end{align*}
We define
\[
U:=\{z \in \C : \langle \Re(z) \rangle < |\Im(z)| < 2 \langle \Re(z) \rangle \}, \qquad V:=\{z \in \C : 0 \leq |\Im(z)| < 2 \langle \Re(z) \rangle \},
\]
and note that the support of $\sigma$ is contained in $V$, while the support of $\sigma_x$ and $\sigma_y$ are contained in $U$. More precisely,
\[
|\sigma_x(z) + i \sigma_y(z)| \lesssim \frac{1}{\langle \Re(z) \rangle} \chi_{U}(z).
\]
Therefore we have the estimate~\cite[Lemma~1]{Davies1995a}\begin{align*}
    \int_{\mathbb{C}} \bigg| &\frac{\partial \tilde{f}}{\partial \overline{z}}\bigg| |z-\lambda_0|^{-1} \cdots |z-\lambda_n|^{-1}dz\\
    & \lesssim \sum_{r=0}^N \int_U |f^{(r)}(\Re(z))| |\Im(z)|^{r-n-1} \langle \Re(z) \rangle^{-1} dz + \int_V |f^{(N+1)}(\Re(z))| |\Im(z)|^{N-n-1} dz\\
    &\lesssim \sum_{r=0}^{N+1} \int_{\R} |f^{(r)}(x)|\langle x \rangle^{r-n-1}dx = \sum_{r=0}^{N+1} \|f\|_{T^{n}(\R), r},
\end{align*}
where the last estimate (integration over the imaginary direction) is justified when $N \geq n+1$.

Hence, if $\|f\|_{T^n(\R),r}<\infty$ for $r=0,\ldots, n+2$, the integral
\begin{align*} 
    \frac{(-1)^{n}}{\pi} \int_{\mathbb{C}} \frac{\partial \tilde{f}}{\partial \overline{z}} (z-\lambda_0)^{-1} \cdots (z-\lambda_n)^{-1}dz
\end{align*}
converges. Now, the identity
\begin{align*}
g^{[n]}(\lambda_0, \ldots, \lambda_n) = \frac{(-1)^{n}}{\pi} \int_{\mathbb{C}} \frac{\partial \tilde{g}}{\partial \overline{z}} (z-\lambda_0)^{-1} \cdots (z-\lambda_n)^{-1}dz, \quad g\in C^\infty_c(\R),
\end{align*}
which holds for compactly supported smooth functions, extends due to the Lebesgue dominated convergence theorem to the functions $f\in C^{n+2}(\R)$ for which $\|f\|_{T^n(\R),r}<\infty$ for $r=0,\ldots, n+2$,
\[
f^{[n]}(\lambda_0, \ldots, \lambda_n) = \frac{(-1)^{n}}{\pi} \int_{\mathbb{C}} \frac{\partial \tilde{f}}{\partial \overline{z}} (z-\lambda_0)^{-1} \cdots (z-\lambda_n)^{-1}dz.
\]
In order to show that this is a decomposition as described in Definition~\ref{def:ProjIntBox}, with $a_j(\lambda_j, z) = (z-\lambda_j)^{-1}$, we will now estimate the expressions
\begin{equation}\label{eq:Sbetanorm}\begin{split}
    \|(z-\cdot)^{-1}\|_{S^{\beta}(\R), k}&=\sup_{\lambda \in \mathbb{R}} \langle \lambda \rangle^{k-\beta} \bigg|\frac{\partial^{k}}{\partial \lambda^{k}} (z-\lambda)^{-1}\bigg| \\
    &\lesssim \bigg(\sup_{\lambda \in \mathbb{R}} \langle \lambda \rangle^{-\beta} |z-\lambda|^{-1}\bigg) \cdot \bigg(\sup_{\lambda \in \mathbb{R}} \langle \lambda \rangle^{k} |z-\lambda|^{-k}\bigg).
\end{split}
\end{equation}

Note that, for $\lambda\in \R$ and $z \in \C \setminus \R$,
\begin{align}\begin{split}\label{eq:Stap1}
    \frac{\langle \lambda \rangle}{|z-\lambda|} = \frac{|\lambda \pm i|}{|z-\lambda|} &\leq 1 + \frac{|z\pm i|}{|z-\lambda|}\\
    &\leq 1+ \frac{\langle z \rangle }{|\Im(z)|},
\end{split}
\end{align}
as  $\min(|z+i|,|z-i|)\leq\langle z\rangle$. Next, we estimate
\[
\sup_{\lambda \in \R} \langle \lambda \rangle^{-\beta} |z-\lambda|^{-1}
\]
for $-1 \leq \beta \leq 0$. We estimate the supremum over $\lambda > 1$, $|\lambda| \leq 1$ and $\lambda < -1$ separately. First, for $|\lambda| \leq 1$ we have $1 \leq \langle \lambda \rangle \leq 2^{1/2}$, and so
\[
\sup_{|\lambda|\leq 1} \langle \lambda \rangle^{-\beta} |z-\lambda|^{-1} 
\lesssim \frac{1}{|\Im(z)|}.
\]
For $\lambda > 1$, we have $\langle \lambda \rangle \leq 2^{1/2} \lambda$, so that
\begin{align*}
    \sup_{\lambda > 1} \langle \lambda \rangle^{-\beta} |z-\lambda|^{-1} &\lesssim 
    \sup_{\lambda > 1-\Re(z)}  \frac{(\lambda+\Re(z))^{-\beta}}{\big(\lambda^2+\Im(z)^2\big)^{1/2}}.
\end{align*}
Writing $v = (\lambda, \Im(z)) \in \R^2$, then by using Cauchy--Schwarz we have
\begin{align*}
    \sup_{\lambda > 1-\Re(z)}  \frac{(\lambda+\Re(z))^{-\beta}}{\big(\lambda^2+\Im(z)^2\big)^{1/2}} & = \sup_{\lambda > 1-\Re(x)}  \frac{(v \cdot (1, \frac{\Re(z)}{\Im(z)}))^{-\beta}}{\|v\|}\\
     &\leq \sup_{\lambda > 1-\Re(z)}  \frac{\| (1, \frac{\Re(z)}{\Im(z)})\|^{-\beta}}{\|v\|^{1+\beta}}\\
    &\leq  \frac{|z|^{-\beta}}{ |\Im(z)|}.
\end{align*}

For $\lambda < -1$ we have a similar estimate, and hence combined we have
\begin{equation}\label{eq:Stap2}
    \sup_{\lambda \in \R}\langle \lambda \rangle^{-\beta} |z-\lambda|^{-1} \lesssim  \frac{1}{|\Im(z)|} \max(1, |z|^{-\beta}) \leq \frac{\langle z \rangle^{-\beta}}{|\Im(z)|}.
\end{equation}

Combining~\eqref{eq:Sbetanorm},~\eqref{eq:Stap1} and~\eqref{eq:Stap2} we get an estimate 
\[
\sup_{\lambda \in \mathbb{R}} \langle \lambda \rangle^{k-\beta} \bigg|\frac{\partial^{k}}{\partial \lambda^{k}} (z-\lambda)^{-1}\bigg| \lesssim \frac{\langle z \rangle^{-\beta}}{|\Im(z)|} \bigg(1+\frac{\langle z \rangle}{|\Im(z)|}\bigg)^k.
\]

Let $-1 \leq \beta_0, \ldots, \beta_n \leq 0$. Taking the inequality above and proceeding as before with $N \geq n+1+\sum_{j=0}^n k_j$, we have
\begin{align*}
    \int_{\mathbb{C}} \bigg| \frac{\partial \tilde{f}}{\partial \overline{z}}\bigg|& \left(\sup_{\lambda_0 \in \mathbb{R}} \langle \lambda_0 \rangle^{k_0-\beta_0} \bigg|\frac{\partial^{k_0}}{\partial \lambda_0^{k_0}} (z-\lambda_0)^{-1}\bigg| \right) \cdots \left(\sup_{\lambda_n \in \mathbb{R}} \langle \lambda_n \rangle^{k_n-\beta_n} \bigg|\frac{\partial^{k_n}}{\partial \lambda_n^{k_n}} (z-\lambda_0)^{-1}\bigg| \right)dz \\
    &\lesssim \sum_{r=0}^{N+1}\int_{\mathbb{R}} |f^{(r)}(x)|  \langle x \rangle^{r-n-1-\sum_{j=0}^n \beta_j} dx < \infty.
\end{align*}
This converges in particular for $\sum_{j=0}^n \beta_j = \alpha - n$. Since $-n-1 \leq \sum_{j=0}^n \beta_j \leq 0$, this choice is possible if $-1 \leq \alpha \leq n$. We have therefore proved for $-1 \leq \alpha \leq n$, and $-1 \leq \beta_0, \ldots, \beta_n \leq 0$ such that $\sum_{j=0}^n \beta_j = \alpha - n$, that 
\begin{align*}
\|f^{[n]}(\lambda_0, \ldots, \lambda_n)&\|_{S^{\beta_0}(\R) \boxtimes_i \cdots \boxtimes_i S^{\beta_n}(\R), k_0, \ldots, k_n} \\
&= \int_{\mathbb{C}} \|  (z-\cdot)^{-1} \|_{S^{\beta_0}(\R), k_0} \cdots \|(z-\cdot)^{-1} \|_{S^{\beta_n}(\R), k_n} \bigg| \frac{\partial \tilde{f}}{\partial \overline{z}}\bigg|dz \\
&\lesssim 
\sum_{r=0}^{n+\sum_{j=0}^n k_j+2} \|f\|_{T^\alpha(\R),r}.
\end{align*}

\item For $f \in T^{\alpha}(\R)$, $-1 \leq \alpha \leq n$, we have by the first part of the lemma that for each $n \in \N$,
\[
f^{[n]} \in S^{\beta_0}(\R) \boxtimes_i \cdots \boxtimes_i S^{\beta_n}(\R)
\]
where each $\beta_j$ can be chosen to lie in the interval $[-1, 0]$, and $\sum \beta_j = \alpha - n$. 

For $f \in T^\alpha(\R)$ with $\alpha \leq -1$, we can write $f  = g \cdot (x+i)^{-k}$ where $g \in T^\beta(\R)$, $-1 \leq \beta \leq 0$ and $k \in \N$.
The Leibniz rule for divided differences dictates
\[
f^{[n]}(\lambda_0, \ldots, \lambda_n) = \sum_{l=0}^{n} g^{[l]}(\lambda_0, \ldots, \lambda_l) \big((x+i)^{-k}\big)^{[n-l]}(\lambda_l, \ldots, \lambda_n).
\]
From part 1 and the explicit form of the divided differences of $\big((x+i)^{-k}\big)^{[n]}$ we therefore conclude that each term is an element of
\[
\sum_{\substack{\beta_0, \ldots, \beta_n \leq 0 \\ \sum \beta_j= \alpha -k-n}}S^{\beta_0}(\R) \boxtimes_i \cdots \boxtimes_i S^{\beta_n}(\R),
\]
with the required estimate of norms. By the same reasoning, this argument goes through for functions $f\in C^{n+\sum_jk_j+2}(\R)$.

\item This follows analogously to assertion 2, by analysing $\big(g(x+i)^k\big)^{[n]}$ for $g \in T^{\alpha}(\R)$ with \\ $-1 \leq \alpha \leq 0$.
\end{enumerate}
\end{proof}

\begin{rem}\label{rem:fToLinfty}
    In particular, Lemma~\ref{L:DivDifT} shows that given any spectral measures $E_0, \ldots, E_n$, we have
    \[
    \| f^{[n]} \|_{L_\infty^{\beta_0}(E_0) \hat{\otimes}_i \cdots \hat{\otimes}_i L_\infty^{\beta_n}(E_n)} \leq \|f^{[n]}\|_{S^{\beta_0}(\R) \boxtimes_i \cdots \boxtimes_i S^{\beta_n}(\R), 0, \ldots, 0} \leq \sum_{k=0}^{n+2} \| f\|_{T^\beta(\R), k} < \infty.
    \]
    For $n = 0$, the space of functions that satisfy this condition closely resembles the space $\mathfrak{F}_m(\R)$ used in~\cite{CareyGesztesy2016} in the context of double operator integrals.
\end{rem}
Lemma~\ref{L:DivDifT} and Remark~\ref{rem:fToLinfty}, combined with Theorem~\ref{T:MainMOIConstruction}, finish the proof of Theorem~\ref{T:MainThmForReal}.

\section{MOI Identities and Asymptotic Expansions}
\label{S:MOIIdentsAndAsympt}
\subsection{MOI Identities and Applications}\label{SS:Identities}
The most important identities for our applications of our multiple operator integrals are the following. These are higher-order analogues of~\eqref{eq:lowneridentity} and~\eqref{eq:commutator}, seeing that $f(H) = T_{f^{[0]}}^H()$.

\begin{prop}\label{P:UMOIcom}
Let $a, X_1,\ldots,X_n \in\op(\Theta)$, let $H_i\in \op^{h_i}(\Theta)$, $h_i > 0$ be symmetric and $\Theta$-elliptic and let $f \in T^\beta(\R)$, $\beta \in \R$. Then
\begin{align}
T^{H_0, \ldots, H_n}_{f^{[n]}}(X_1,\ldots,X_j,aX_{j+1},\ldots,X_n)-T_{f^{[n]}}^{H_0, \ldots, H_n}&(X_1,\ldots,X_ja,X_{j+1},\ldots,X_n)\label{eq:comeq1}\\
&=T_{f^{[n+1]}}^{H_0, \ldots, H_j, H_j, \ldots, H_n}(X_1,\ldots,X_j,[H_j,a],X_{j+1},\ldots,X_n);\notag\\
    T^{H_0, \ldots, H_n}_{f^{[n]}}(aX_1,\ldots,X_n) - a T^{H_0, \ldots, H_n}_{f^{[n]}}(X_1,\ldots,X_n) &= T^{H_0, H_0, H_1, \ldots, H_n}_{f^{[n+1]}}([H_0,a],X_1,\ldots,X_n) ;\\
    T^{H_0, \ldots, H_n}_{f^{[n]}}(X_1,\ldots,X_n)a -  T^{H_0, \ldots, H_n}_{f^{[n]}}(X_1,\ldots,X_na) &= T^{H_0, \ldots, H_n, H_n}_{f^{[n+1]}}(X_1,\ldots,X_n,[H_n,a]).
\end{align}
Moreover, for $A \in \op^a(\Theta), a>0$, $B \in \op^b(\Theta), b>0$ symmetric and $\Theta$-elliptic,
\begin{align}\label{eq:SuperscriptDiff}
    \begin{split}
        &T^{H_0, \ldots, H_{j-1}, A, H_{j+1}, \ldots, H_n}_{f^{[n]}}(X_1, \ldots, X_n) - T^{H_0, \ldots, H_{j-1}, B, H_{j+1}, \ldots, H_n}_{f^{[n]}}(X_1, \ldots, X_n)\\
        &=T^{H_0, \ldots, H_{j-1}, A, B, H_{j+1}, \ldots, H_n}_{f^{[n+1]}}(X_1, \ldots, X_j, A-B, X_{j+1}, \ldots, X_n).
    \end{split}
\end{align}
\end{prop}
\begin{proof}
Note that for $f\in T^\beta(\R)$, the multiple operator integrals appearing above are well-defined through Theorem~\ref{T:MainThmForReal}.

We prove equation~\ref{eq:comeq1}, the others follow analogously. Write
\begin{align*}
    F_j(\lambda_0, \ldots, \lambda_{n+1}) &:= f^{[n]}(\lambda_0, \ldots, \lambda_{j-1}, \lambda_{j+1}, \ldots \lambda_{n+1});\\
    F_{j+1}(\lambda_0, \ldots, \lambda_{n+1}) &:= f^{[n]}(\lambda_0, \ldots, \lambda_{j}, \lambda_{j+2}, \ldots \lambda_{n+1}).
\end{align*}
Observe that
\[
F_{j+1}(\lambda_0, \ldots, \lambda_{n+1}) - F_j(\lambda_0, \ldots, \lambda_{n+1}) = (\lambda_j - \lambda_{j+1}) f^{[n+1]}(\lambda_0, \ldots, \lambda_{n+1}).
\]
Hence, 
\begin{align*}
    T_{f^{[n+1]}}^{H_0, \ldots, H_j, H_j, \ldots, H_n}&(X_1,\ldots,X_j,[H_j,a],X_{j+1},\ldots,X_n)\\
    &= T_{(\lambda_j - \lambda_{j+1}) f^{[n+1]}}^{H_0, \ldots, H_j, H_j, \ldots, H_n}(X_1,\ldots,X_j,a,X_{j+1},\ldots,X_n)\\
    &= T_{F_{j+1}}^{H_0, \ldots, H_j, H_j, \ldots, H_n}(X_1,\ldots,X_j,a,X_{j+1},\ldots,X_n) \\
    &\quad - T_{F_j}^{H_0, \ldots, H_j, H_j, \ldots, H_n}(X_1,\ldots,X_j,a,X_{j+1},\ldots,X_n) \\
    &= T^{H_0, \ldots, H_n}_{f^{[n]}}(X_1,\ldots,X_j,aX_{j+1},\ldots,X_n)-T_{f^{[n]}}^{H_0, \ldots, H_n}(X_1,\ldots,X_ja,X_{j+1},\ldots,X_n).
\end{align*}
\end{proof}

With these identities in hand, we can show that the MOI constructed in the previous section is an element of $\OP(\Theta)$ if all its components are and the symbol is a divided difference.
\begin{proof}[Proof of Theorem~\ref{T:MOIpsi}]
     Let $H_i \in \OP^{h}(\Theta)$, $h>0$, be symmetric and $\Theta$-elliptic and take $f \in T^\beta(\R)$, $\beta \in \R$. For operators $X_i \in \OP^{r_i}(\Theta)$, $r:= \sum_{i=1}^n r_i$, we have that
     \[
     T_{f^{[n]}}^{H_0, \ldots, H_n}(X_1, \ldots, X_n) \in \op^{h(\beta-n)+ r}(\Theta)
     \]
     due to Lemma~\ref{L:DivDifT}.
     Now, taking $n=1$ to ease notation, using Proposition~\ref{P:UMOIcom} gives that
     \[
     [\Theta, T^{H_0, H_1}_{f^{[1]}}(X_1)] = T^{H_0, H_0, H_1}_{f^{[2]}}([\Theta, H_0],X_1) + T^{H_0, H_1}_{f^{[1]}}([\Theta,X_1]) + T^{H_0, H_1, H_1}_{f^{[2]}}(X_1, [\Theta, H_1]).
     \]
     As $[\Theta, H_i]\in \op^h(\Theta)$ and $[\Theta, X_1] \in \op^{r}(\Theta)$, Lemma~\ref{L:DivDifT} combined with Theorem~\ref{T:MainMOIConstruction} gives that
     \[
     T^{H_0, H_0, H_1}_{f^{[2]}}([\Theta, H_0],X_1), T^{H_0,H_1}_{f^{[1]}}([\Theta,X_1]), T^{H_0, H_1, H_1}_{f^{[2]}}(X_1, [\Theta, H_1]) \in \op^{h(\beta-n)+r}(\Theta).
     \]
     Higher commutators and $n > 1$ follow analogously. Hence
     \[
     T_{f^{[n]}}^{H_0,\ldots, H_n}(X_1, \ldots, X_n) \in \OP^{h(\beta-n)+ r}(\Theta).
     \]
\end{proof}

In the setting that $\Theta^{-1} \in \mathcal{L}_s$, $s>0$, it is immediate from Theorem~\ref{T:MainMOIConstruction} and Theorem~\ref{T:MainFunctCalc} that for $H_0, \ldots, H_n \in \op^{h}(\Theta), h >0$ symmetric and $\Theta$-elliptic, and $f\in T^\beta(\R)$, the multiple operator integral
\begin{equation}\label{eq:TypicalMOI}
    T_{f^{[n]}}^{H_0, \ldots, H_n}(X_1, \ldots, X_n) \in \op^{(\beta-n)h + r}(\Theta)
\end{equation}
will be a trace-class operator when considered as an operator on $\H$ if $\beta$ is small enough. Namely, we have
\[
\|A\|_{1} \leq \|\Theta^{-s}\|_1 \|A\|_{\H^0 \to \H^{s}}.
\] 
For asymptotic expansions of trace formulas, it will be useful to make a more detailed analysis.

\begin{prop}\label{P:TraceEstimatewitht}
Let $H_i \in \op^h(\Theta)$ for a fixed $h > 0$ be symmetric and $\Theta$-elliptic operators. If $X_i \in \op^{r_i}(\Theta)$, $r:= \sum_{i=1}^n r_i$, $f \in T^\alpha(\mathbb{R})$ with $\alpha \leq n$, $t\leq 1$, then
    \[
    \| T^{tH_0, \ldots, tH_n}_{f^{[n]}}(X_1, \ldots, X_n) \|_{\H^{q+r+(\alpha-n)h} \to \H^q} \lesssim t^{\alpha-n}, \quad q\in \R.
    \]
Rephrased, if $f \in T^{\frac{u-r}{h}+n}(\mathbb{R})$, $u\leq r$, then
\[
    \| T^{tH_0, \ldots, tH_n}_{f^{[n]}}(X_1, \ldots, X_n) \|_{\H^{q+u} \to \H^q} \lesssim t^{\frac{u-r}{h}}, \quad q\in \R.
    \]
\end{prop}
\begin{proof}
Lemma~\ref{L:DivDifT} gives that
\[
      f^{[n]}\in \sum_{\beta_0 + \cdots + \beta_n = \alpha-n} S^{\beta_0}(\R) \boxtimes_i \cdots \boxtimes_i S^{\beta_n}(\R),
\]
where each $\beta_i \leq 0$ since $\alpha \leq n$.
   Consider one of the summands, $\phi \in S^{\beta_0}(\R) \boxtimes_i \cdots \boxtimes_i S^{\beta_n}(\R)$. Then 
   \[
   \phi(\lambda_0, \ldots, \lambda_n) = (\lambda_0 + i)^{\beta_0} \cdots (\lambda_n + i)^{\beta_n} \cdot \psi(\lambda_0, \ldots, \lambda_n),
   \]
   for a function $\psi \in S^0(\R) \boxtimes_i \cdots \boxtimes_i S^0(\R)$, and thus, by Remark~\ref{rem:SipvL} and Theorem~\ref{T:MainMOIConstruction},
\begin{align*}
    T^{tH_0, \ldots, tH_n}_{\phi}&(X_1, \ldots, X_n)\\
    &=  T^{tH_0, \ldots, tH_n}_{\psi}((tH_0+i)^{\beta_0} X_1 (tH_1+i)^{\beta_1}, X_2(tH_2+i)^{\beta_2},\ldots, X_n (tH_n+i)^{\beta_n}).
\end{align*}
Corollary~\ref{C:IndepNorm} and Theorem~\ref{T:MainMOIConstruction} give that
\begin{align*}
    &\|  T^{tH_0, \ldots, tH_n}_{\psi}((tH_0+i)^{\beta_0}  X_1 (tH_1+i)^{\beta_1}, X_2(tH_2+i)^{\beta_2},\ldots, X_n (tH_n+i)^{\beta_n}) \|_{\H^{q+r+(\alpha-n)h}\to\H^q} \\
    &\lesssim\| (tH_0+i)^{\beta_0}\|_{\H^{q_0 + \beta_0h} \to \H^{q_0}} \cdots \| (tH_n+i)^{\beta_n}\|_{\H^{q_n+\beta_nh}\to \H^{q_n}},
\end{align*}
for $q_i$ some real numbers. Theorem~\ref{T:MainFunctCalc} gives that 
\[
\|(tH_j+i)^{\beta_j}\|_{\H^{q_j+\beta_jh}\to \H^{q_j}} \lesssim \sup_{x\in \R} |(tx+i)^{\beta_j} | \langle x\rangle^{-\beta_j} \lesssim t^{\beta_j}.
\]
Therefore, 
\begin{align*}
    \|  T^{tH_0, \ldots, tH_n}_{\phi}(X_1, X_2,\ldots, X_n )& \|_{\H^{q+r+(\alpha-n)h} \to \H^q} \lesssim t^{\beta_0 + \cdots + \beta_k}= t^{\alpha-n}.
\end{align*}
\end{proof}

\subsection{Asymptotic Expansions}
\label{SS:AsympExp}
Through the identities proved in the previous section, the theory of multiple operator integrals lends itself well for establishing asymptotic expansions of operators. As an immediate example, we first prove a Taylor expansion for pseudodifferential operators, which is the first part of Theorem~\ref{T:AsypExpansions}.
\begin{thm}\label{T:Taylor}
     Let $f \in T^\beta(\R)$, $H \in \op^h(\Theta)$, $h>0$ $\Theta$-elliptic and symmetric, and let $V \in \op^r(\Theta)$ be symmetric. If the order of the perturbation $V$ is strictly smaller than that of $H$, i.e. $r < h$, we have
    \begin{align}\label{eq:Taylor asymptotic formula}
    f(H + V) \sim &\sum_{n=0}^\infty T^{H}_{f^{[n]}}(V, \ldots, V),
    \end{align}
    in the sense that
    \begin{align*}
        &f(H+V) - \sum_{n=0}^N T^{H}_{f^{[n]}}(V, \ldots, V)\in \op^{m_N}(\Theta)
    \end{align*}
    with $m_N \downarrow -\infty$.
\end{thm}
\begin{proof}
Using the last part of Proposition~\ref{P:UMOIcom} with $A = H+V$, $B = H$, we have
\begin{align*}
    f(H+V) - f(H) &= T^{H+V}_{f^{[0]}}() - T^H_{f^{[0]}}()\\
    &= T_{f^{[1]}}^{H+V, H}(V).
\end{align*}
Repeating the argument, we get
    \begin{align*}
            f(H+V) - &\sum_{n=0}^N  T^{H}_{f^{[n]}}(V, \ldots, V) = T_{f^{[N+1]}}^{H+V, H, \ldots, H}(V, \ldots, V).
    \end{align*}
Now, if $h > r$, Theorem~\ref{T:MainThmForReal} gives that
\[
T_{f^{[N+1]}}^{H+V, H, \ldots, H}(V, \ldots, V) \in \op^{(\beta-N-1)h+Nr}(\Theta),
\]
with 
\[
(\beta-N-1)h+Nr = N(r-h) + (\beta-1)h \downarrow -\infty.
\]
\end{proof}

Note that, if $H$ and $V$ are commuting numbers, \eqref{eq:Taylor asymptotic formula} recovers the classic Taylor expansion formula $f(H+V)\sim\sum_{n=0}^\infty \frac{f^{(n)}(H)}{n!}V^n$. The noncommutative Taylor expansion \eqref{eq:Taylor asymptotic formula} features prominently in \cite{Hansen,vNvS21a,Suijlekom2011,Skripka2018,lesch2017divided} for bounded operators $V$.

Each term in the Taylor expansion can itself be expanded as follows. Recall that we write $\delta_H(X) := [H,X]$, $\delta_H^n(X) := \delta_H(\cdots \delta_H(\delta_H(X)) \cdots)$. The following proposition and corollary finish the proof of Theorem~\ref{T:AsypExpansions}.

\begin{prop}\label{P:Expansion} Let $X_i \in \op^{r_i}(\Theta)$, $H\in \op^h(\Theta)$, $h >0$ symmetric and $\Theta$-elliptic, and $f\in T^\beta(\R)$. Then
    \begin{align*}
          T^{H}_{f^{[n]}}(X_1, \dots, X_n) &= \sum_{m=0}^{N}  \sum_{m_1 + \dots + m_n = m} \frac{C_{m_1, \dots, m_n}}{(n+m)!} \delta_H^{m_1}(X_1) \cdots \delta_H^{m_n}(X_n) f^{(n+m)}(H) \\
        &\quad +   S^n_{N}(X_1, \dots, X_n),
    \end{align*}
    where 
    \[
    C_{m_1, \dots, m_n} :=  \prod_{j=1}^{n} \binom{j+m_1 + \dots + m_{j}-1}{m_j}
    \]
    and the remainder $   S^n_{N}(X_1, \dots, X_n)$ is a sum of terms of the form
    \[
 \delta_H^{m_1}(X_1)\cdots \delta_H^{m_k}(X_{k})   T^{H}_{f^{[n+N+1]}}(1, \ldots, 1, \delta_H^{N+1-m_1 - \cdots - m_k}(X_{k+1}), 1, \ldots, 1, X_{k+2}, \ldots, X_n).
\]
If the commutators $\delta_H^k(X_j)$ have a lower order than the expected $r_j+kh$, explicitly if
    \[
    \delta_H^k(X_j) \in \op^{r_j+k(h-\varepsilon)}(\Theta)
    \]
    for some $\varepsilon > 0$, then the above gives an asymptotic expansion
    \[
        T^{H}_{f^{[n]}}(X_1, \dots, X_n) \sim \sum_{m=0}^{\infty}  \sum_{m_1 + \dots + m_n = m} \frac{C_{m_1, \dots, m_n}}{(n+m)!} \delta_H^{m_1}(X_1) \cdots \delta_H^{m_n}(X_n) f^{(n+m)}(H),
    \]
    in the sense that the remainder term 
    \[
    S^n_{N}(X_1, \ldots, X_n) \in \op^{m_N}(\Theta)
    \]
    with $m_N = \sum_j r_j + (\beta-n)h-\varepsilon (N+1) \downarrow -\infty$.
\end{prop}
\begin{proof}
As the proof is a lengthy combinatorial exercise, which consists of repeatedly applying Proposition~\ref{P:UMOIcom}, we have elected to move the proof of this proposition to Appendix~\ref{S:AppCombinComp}.
\end{proof}

\begin{cor}\label{C:CombiExpansion}
    Under the combined assumptions of Theorem~\ref{T:Taylor} and Proposition~\ref{P:Expansion}, we get
    \[
    f(H+V) \sim  \sum_{n,m=0}^\infty \sum_{m_1 + \dots + m_n = m} \frac{C_{m_1, \dots, m_n}}{(n+m)!} \delta_H^{m_1}(V) \cdots \delta_H^{m_n}(V) f^{(n+m)}(H).
    \]
\end{cor}
\begin{proof}
    For $N, M > 0$, Theorem~\ref{T:Taylor} and Proposition~\ref{P:Expansion} give
    \begin{align*}
        f(H+V) &= \sum_{n=0}^{N} \sum_{m=0}^{M} \sum_{m_1 + \dots + m_n = m} \frac{C_{m_1, \dots, m_n}}{(n+m)!} \delta_H^{m_1}(V) \cdots \delta_H^{m_n}(V) f^{(n+m)}(H)\\
        &\quad + \op^{N(r-h)+(\beta-1) h}(\Theta) + \sum_{n=0}^N \op^{nr + (\beta-n) h-\varepsilon(M+1)}(\Theta).
    \end{align*}
    As $(N,M) \to (\infty, \infty)$, we see that the order of the remainder decreases to $-\infty$.
\end{proof}

A version of Proposition~\ref{P:Expansion} for classical pseudodifferential operators (i.e., on a manifold) appears as \cite[Theorem 1]{Paycha2007}. It should also be noted that the combinatorial or algebraic manipulations needed to prove Proposition~\ref{P:Expansion} form a part of every proof of the local index formula, see for example~\cite[Lemma~6.11]{CPRS1}, \cite[Equation~(71)]{ConnesMoscovici1995}, \cite[Lemma~2.12]{Higson2003}. Namely, the cocycle that lies at the heart of the local index formula is the residue of the trace of an expansion of the above kind. 

This observation leads us naturally to the topic of asymptotic expansions of trace formulas. In various contexts of noncommutative geometry and beyond, expansions are studied of the kind
\begin{equation}\label{eq:AsympExp1}
    \Tr(f(tH+tV)) \underset{t\downarrow 0}{\sim} \sum_{k=0}^\infty c_k t^{r_k},
\end{equation}
for an increasing sequence $r_k \uparrow \infty$ and constants $c_k \in \mathbb{C}$, which means that as $t \downarrow 0$
\[
\Tr (f(tH+tV)) = \sum_{k=0}^N c_k t^{r_k} + O(t^{r_{N+1}})
\]
for every $N \in \mathbb{Z}_{\geq 0}$. Or, more generally (c.f.~\cite{EcksteinIochum2018}),
\begin{equation}\label{eq:AsympExp2}
     \Tr(f(tH+tV)) \underset{t\downarrow 0}{\sim} \sum_{k=0}^\infty \rho_k(t),
\end{equation}
where $\rho_k(t) = O(t^{r_k})$ and the asymptotic expansion means that for every $N\in \mathbb{Z}_{\geq 0}$
\[
\Tr(f(tH+tV)) = \sum_{k=0}^N \rho_k(t) + O(t^{r_{N+1}}).
\]

We will target a specific open problem posed in some form by Eckstein and Iochum in~\cite{EcksteinIochum2018}. 
Given a spectral triple $(\A, \H, D)$ it is a common assumption to require the existence of an asymptotic expansion as $t \downarrow 0$ of
\[
\Tr(ae^{-tD^2}),
\]
where $a \in \A$. It is currently not clear 
whether the existence of asymptotic expansions of
\[
\Tr(ae^{-t(D+V)^2})
\]
can be deduced for suitable perturbations $V$ from this, and whether it could be enough for $\Tr(e^{-tD^2})$ to admit an asymptotic expansion.

\begin{thm}{\cite[Theorem~3.2]{EcksteinIochum2018}}\label{T:ExpansionToZeta}
    For a bounded operator $a$ and invertible positive operator $D$ such that $D^{-1} \in \mathcal{L}_s$, $s >0$, the existence of an asymptotic expansion 
    \[
    \Tr(ae^{-tD^2}) \underset{t \downarrow 0}{\sim} \sum_{k=0}^\infty \rho_k(t),
    \]
    where 
    \[
    \rho_k(t) := \sum_{z \in X_k} \bigg( \sum_{n=0}^d c_{n,k} \log^n t \bigg) t^{-z},
    \]
    with $c_{n,k} \in \C$ and for suitable sets $X_k \subset \C$ (for details, see~\cite[Theorem~3.2]{EcksteinIochum2018}),
    implies the existence of a meromorphic continuation of
    \[
    \zeta_{D^2,a}(s) := \Tr(a |D|^{-2s})
    \]
    to the complex plane, with poles of order at most $d+1$ located at points in $\bigcup_{k=0}^\infty X_k \subset \C$.
\end{thm}

We first show that in general the asymptotic expansion of $\Tr(e^{-tD^2})$ provides no control over the expansions of $\Tr(a e^{-tD^2})$.
\begin{ex}\label{E:Counterexample}
    Let $\A = \ell_{\infty}(\Z_{\geq 1}),$ $\H = \ell_2(\Z_{\geq 1})$ where $\A$ is represented
    on $\H$ by pointwise multiplication, and let $D$ be the diagonal operator on $\H$ given by
    \[
        De_n = ne_n,\quad n\geq 1.
    \]
This is a spectral triple for trivial reasons: $\A$ acts on $\H$ by bounded operators,
and $[D,a]=0$ for all $a\in\A.$
Despite being atypical, $(\A,\H,D)$ satisfies most of the assumptions commonly made in the literature in terms of smoothness or summability.
The algebra $\A$ is not separable, but all of the following arguments can be performed in a separable (even finite dimensional) subalgebra of $\A.$

It is a classical result that we have the asymptotic expansion
\[
\Tr(e^{-tD^2})= \sum_{n=1}^\infty e^{-tn^2} \underset{t \downarrow 0}{\sim} \frac{\sqrt{\pi}}{2} t^{-\frac{1}{2}} - \frac{1}{2},
\]
see for example~\cite[Lemma~3.1.3]{Gilkey2004}. 
Nonetheless, the functions $\Tr(a|D|^{-2s})$ for $a\in \A$ are very badly behaved. For example, let
\[
    a:= \sum_{n=2}^\infty \frac{1}{\log n} e_n \in \A,
\]
so that
\[
\zeta_{a,D^2}(s)= \Tr(aD^{-2s}) = \sum_{n=2}^\infty  \frac{1}{\log n}n^{-2s},\quad \Re(s)>\frac{1}{2},
\]
which is holomorphic on $\Re(s) > \frac{1}{2}$. Now,
\begin{align*}
    \frac{d}{ds}\zeta_{a,D^2}(s) &= -2\sum_{n=2}^\infty n^{-2s}  \\
    &=2-2\zeta(2s),\quad \Re(s)>\frac{1}{2},
\end{align*}
where $\zeta$ is the Riemann zeta function which has a simple pole at $1$. Therefore, $\zeta_{a,D^2}(s) + \log(2s-1)$ is the antiderivative of an entire function, which implies that $\zeta_{a,D^2}(s) = -\log(2s-1) + f(s)$ where $f(s)$ is entire~\cite[Theorem~10.14]{PapaRudin}.
We conclude that $\zeta_{a,D^2}$ does not admit a meromorphic extension to the complex plane, and thus $\Tr(ae^{-tD^2})$ does not have an asymptotic expansion as $t\downarrow 0$ of the type in Theorem~\ref{T:ExpansionToZeta}.

An even more pathological example is 
\[
    b := \sum_{n=2}^\infty \frac{\Lambda(n)}{\log(n)}e_n
\]
where $\Lambda$ is the von Mangoldt function which satisfies
\[
\Lambda(n) :=
\begin{cases}
    \log(p) &\text{ if } n=p^k \text{ for } p \text{ prime};\\
    0 &\text{ otherwise}.
\end{cases}
\] 
A classical formula asserts that~\cite[p.4]{Titchmarsh1986}
\[
    \zeta_{b,D}(s)=\Tr(b|D|^{-s}) = \sum_{n=2}^\infty \frac{\Lambda(n)}{\log(n)}n^{-s} = \log \zeta(s), \quad \Re(s) > 1
\]
which is badly behaved at every zero of $\zeta$ and at $s= 1$.
\end{ex}

To study the asymptotic expansion of expressions like
\[
\Tr(a e^{-t(D+V)^2}),
\]
we will use a modified version of Theorem~\ref{T:Taylor} and Proposition~\ref{P:Expansion}. For this purpose, we first analyse the remainder in the Taylor expansion in Theorem~\ref{T:Taylor} more precisely.

\begin{prop}\label{P:TaylorTrace}Let $\Theta^{-1} \in \mathcal{L}_{s}$, $s > 0$, $f \in T^\beta(\R)$, $H \in \op^h(\Theta), h> 0$ $\Theta$-elliptic and symmetric and $V \in \op^r(\Theta)$ symmetric. Let $h > r \geq 0$ and $\beta < -\frac{s}{h}$.
    For every $N \in \N$, we have as $t \downarrow 0$,
    \[
    \Tr(f(tH+tV)) = \sum_{n=0}^{N} t^{n} \Tr( T^{tH}_{f^{[n]}}(V, \ldots, V) ) + O(t^{(N+1)(1-\frac{r}{h}) - \frac{s}{h}}).
    \]
\end{prop}
\begin{proof}
    The proof of Theorem~\ref{T:Taylor} gives that
    \[
    f(tH + tV) = \sum_{n=0}^{N} T_{f^{[n]}}^{tH}(tV, \ldots, tV) + T_{f^{[N+1]}}^{tH + tV, tH, \ldots, tH}(tV, \ldots, tV).
    \]
    The condition $\beta < -\frac{s}{h}$ assures that all terms on the left and right-hand side are trace-class (cf. \eqref{eq:TypicalMOI}). Furthermore, we have $f \in T^\beta(\R)\subseteq T^{(N+1)(1-\frac{r}{h})-\frac{s}{h}}(\R)$ so that Proposition~\ref{P:TraceEstimatewitht} provides that
    \begin{align*}
        \|T_{f^{[N+1]}}^{tH + tV, tH, \ldots, tH}&(tV, \ldots, tV)\|_1 \\
        &\lesssim t^{N+1} \|T_{f^{[N+1]}}^{tH + tV, tH, \ldots, tH}(V, \ldots, V)\|_{\H^{-s}\to\H^0}\\
        &\lesssim t^{(N+1)(1-\frac{r}{h}) - \frac{s}{h}}.
    \end{align*}
\end{proof}

This proposition makes it clear that to determine the coefficients of asymptotic expansions of the type~\eqref{eq:AsympExp1} or~\eqref{eq:AsympExp2}, it suffices to study the asymptotic expansions of the multiple operator integral
\[
\Tr( T^{tH}_{f^{[n]}}(V, \ldots, V) ),
\]
which we do with Proposition~\ref{P:Expansion}.

\begin{prop}\label{P:GeneralAsymp}
    Let $\Theta^{-1} \in \mathcal{L}_{s}$, $s > 0$, $f \in T^\beta(\R)$, $H\in \op^h(\Theta)$, $h>0$ symmetric and $\Theta$-elliptic, $V \in \op^r(\Theta)$ symmetric. If $h > r \geq0 $, $\beta \leq -\frac{s}{h}$, and $\delta_H^n(V) \in \op^{r+n(h-\varepsilon)}(\Theta)$, then $\Tr(f(tH+tV)) $ admits an asymptotic expansion as $t\downarrow 0$ of type~\eqref{eq:AsympExp2} given by
    \[
    \Tr(f(tH+tV)) = \sum_{n=0}^{N} \sum_{m=0}^{N}   \sum_{m_1 + \dots + m_n = m} t^{n+m}\frac{C_{m_1, \dots, m_n}}{(n+m)!} \Tr\big( \delta_H^{m_1}(V) \cdots \delta_H^{m_n}(V) f^{(n+m)}(tH)\big) + O(t^{m_N}),
    \]
    where $m_N := (N+1)\min\big(\frac{\varepsilon}{h}, (1-\frac{r}{h}) \big) - \frac{s}{h},$ so that $m_N \uparrow \infty$ as $N \to \infty$.
\end{prop}
\begin{proof}
    Combining Propositions~\ref{P:TaylorTrace} and~\ref{P:Expansion}, we have that
    \begin{align*}
        \Tr(f(tH+tV)) &= \sum_{n=0}^{N} \sum_{m=0}^{N}  t^{n+m} \sum_{m_1 + \dots + m_n = m} \frac{C_{m_1, \dots, m_n}}{(n+m)!} \Tr\big( \delta_H^{m_1}(V) \cdots \delta_H^{m_n}(V) f^{(n+m)}(tH)\big)\\
        &\quad + \sum_{n=0}^{N} t^n \Tr(S^n_{N,t}(V, \ldots, V))+ O(t^{(N+1)(1-\frac{r}{h}) - \frac{s}{h}}),
    \end{align*}
    where $S^n_{N,t}(V, \ldots, V)$ is a sum of terms of the form
    \[
    t^{N+1}   \delta_H^{m_1}(V)\cdots \delta_H^{m_k}(V)   T^{tH}_{f^{[n+N+1]}}(1, \ldots, 1, \delta_H^{N+1-m_1 - \cdots - m_k}(V), 1, \ldots, 1, V, \ldots, V).
    \]
    We then estimate
    \begin{align*}
        &\big \|t^{N+1}   \delta_H^{m_1}(V)\cdots \delta_H^{m_k}(V)   T^{tH}_{f^{[n+N+1]}}(1, \ldots, 1, \delta_H^{N+1-m_1 - \cdots - m_k}(V), 1, \ldots, 1, V, \ldots, V)\big\|_1\\
        &\lesssim t^{N+1} \big \|  \delta_H^{m_1}(V)\cdots \delta_H^{m_k}(V)   T^{tH}_{f^{[n+N+1]}}(1, \ldots, 1, \delta_H^{N+1-m_1 - \cdots - m_k}(V), 1, \ldots, 1, V, \ldots, V)\big\|_{\H^{-s}\to\H^0}\\
        &\leq t^{N+1} \big\|\delta_H^{m_1}(V)\cdots \delta_H^{m_k}(V) \big\|_{\H^{kr+(m_1 + \cdots + m_k)(h-\varepsilon)} \to \H^{0}}\\
        &\quad \times \big\| T^{tH}_{f^{[n+N+1]}}(1, \ldots, 1, \delta_H^{N+1-m_1 - \cdots - m_k}(V), 1, \ldots, 1, V, \ldots, V)\big\|_{\H^{-s}\to \H^{kr+(m_1 + \cdots + m_k)(h-\varepsilon)}}.
    \end{align*}
    Applying Proposition~\ref{P:TraceEstimatewitht} then provides that
    \begin{align*}
        \| S^n_{N, t}(V, \ldots, V)\|_1 &\lesssim t^{N+1-\frac{s+kr+(m_1 + \cdots + m_k)(h-\varepsilon)}{h}-\frac{(n-k)r+(N+1-m_1 - \cdots - m_k)(h-\varepsilon)}{h}}\\
        &\quad = t^{-\frac{s+nr}{h}+(N+1)\frac{\varepsilon}{h}},
    \end{align*}
    and hence
    \begin{align*}
        \sum_{n=0}^{N} t^n \Tr(S^n_{N,t}(V, \ldots, V)) \lesssim \sum_{n=0}^{N} t^{n(1-\frac{r}{h})+(N+1)\frac{\varepsilon}{h}-\frac{s}{h}}.
    \end{align*}
    Defining
    \[
    m_N := (N+1)\min\bigg(\frac{\varepsilon}{h}, \big(1-\frac{r}{h}\big) \bigg) - \frac{s}{h}
    \]
    concludes the proof.
\end{proof}

\begin{proof}[Proof of Theorem~\ref{T:AsympExpNCG}]
Given a regular $s$-summable spectral triple  $(\A,\H,D)$, write $\Theta := (1+D^2)^{1/2}$. Let $V, P \in \mB$, $V$ self-adjoint and bounded, where $\mB$ is the algebra generated by $\A$ and $D$. If $f \in T^{\beta}(\R)$ with $\beta < -s$, Proposition~\ref{P:GeneralAsymp} immediately gives that
  \[
    \Tr(f(tD+tV)) = \sum_{n=0}^{N} \sum_{m=0}^{N}   \sum_{m_1 + \dots + m_n = m} t^{n+m}\frac{C_{m_1, \dots, m_n}}{(n+m)!} \Tr\big( \delta_D^{m_1}(V) \cdots \delta_D^{m_n}(V) f^{(n+m)}(tD)\big) + O(t^{N+1-s}).
    \]

Regarding the expansion of $\Tr(Pe^{-t(D+V)^2})$ we have that $D^2 \in \OP^2(\Theta)$, $A \in \OP^1(\Theta)$ since $(\A, \H, D)$ is regular. As $[D^2, A] = [\Theta^2, A] = \Theta[\Theta, A] + [\Theta,A] \Theta$, we have that $A^{(m)} \in \OP^{1+m}(\Theta)$. Furthermore, $\mB \subseteq \op(\Theta)$. The proof of this corollary is then the same as the proof of Proposition~\ref{P:GeneralAsymp}. Filling in
\[
m_N = (N+1)\min\bigg(\frac{\varepsilon}{h}, \big(1-\frac{r}{h}\big) \bigg) - \frac{s}{h} = \frac{N+1-s}{2}
\]
gives the order of the error term.

For the expansion of $\Tr(Pe^{-t|D+V|})$, while $|D| \in \OP^1(\Theta)$ due to Theorem~\ref{T:MainFunctCalc}, to conclude something similar for $|D+V|$ we have to do more work. Note that $D$ has discrete spectrum since $(1+D^2)^{-1/2} \in \mathcal{L}_s$. If $V \in \OP^0(\Theta)$ is self-adjoint, then $D+V$ has real discrete spectrum too. Hence we can modify the function $x \mapsto |x|$ slightly on a small neighbourhood around $x=0$ to get a smooth function $f$ which has the property that $f(x) = |x|$ on $\sigma(D+V) \cup \sigma(D)$, and $f \in S^1(\R)$ since the second and higher derivatives of $f$ are all compactly supported. Using Theorem~\ref{T:MOIpsi} and the observation that $S^1(\R) \subseteq T^{1+\varepsilon}(\R)$ for all $\varepsilon > 0$, we have
\[
|D+V| \in \OP^{1+\varepsilon}(\Theta), \ \ |D+V|-|D| = T^{D+V, D}_{f^{[1]}}(V) \in \OP^{\varepsilon}(\Theta).
\]
Therefore, we get
\begin{align*}
    &\Tr(Pe^{-t|D+V|}) \\
    &= \sum_{n = 0}^{N} \sum_{m=0}^{N} \sum_{m_1 + \dots + m_n = m}  (-t)^{n+m} \frac{C_{m_1, \dots, m_n}}{(n+m)!} \Tr(P\delta_{|D|}^{m_1}(B)\cdots \delta_{|D|}^{m_n}(B)\exp(-t|D|)) + O(t^{(N+1)(1-\varepsilon) - s}),
\end{align*}
where $B := |D+V|-|D|$ and $\varepsilon > 0$ can be chosen arbitrarily small.
\end{proof}

Apart from providing a pertubative expansion of the spectral action, Theorem~\ref{T:AsympExpNCG} shows that if for all $P \in \mB$ we have an expansion
\begin{equation}\label{eq:ExistenceAsExp}
\Tr (P e^{-tD^2}) \underset{t\downarrow 0}{\sim} \sum_{k=0}^\infty c_k(P) t^{r_k}
\end{equation}
for constants $c_k(P) \in \C$, then there exist constants $c_k(P,V) \in \C$ such that
\[
\Tr(P e^{-t(D+V)^2}) \underset{t\downarrow 0}{\sim} \sum_{k=0}^\infty c_k(P, V) t^{r_k}.
\]
Similarly, if each $\Tr (P e^{-tD^2})$ admits an asymptotic expansions of the type in Theorem~\ref{T:ExpansionToZeta}, then $\Tr (P e^{-t(D+V)^2})$ does too.

\begin{rem}
    Theorem~\ref{T:AsympExpNCG} can be modified to work for non-unital spectral triples. Given a spectral triple $(\A, \H, D)$ with non-unital algebra $\A$ writing $\Theta = (1+D^2)^{1/2}$, if one assumes instead of $\Theta^{-1} \in \mathcal{L}_s$ that there exists $p \geq 1$ such that $a \Theta^{-s} \in \mathcal{L}_1$ for all $a \in \A \cup [D,\A]$ and $s > p$ as is proposed in~\cite{CGRS2}, then we also have $a \cdot \op^{-s}(\Theta) \in \mathcal{L}_1$ for $s > p$. It follows that as $t\downarrow 0$
    \[
    \Tr(af(tD+tV)) = \sum_{n=0}^{N} \sum_{m=0}^{N}  t^{n+m} \sum_{m_1 + \dots + m_n = m} \frac{C_{m_1, \dots, m_n}}{(n+m)!} \Tr\big(a \delta_D^{m_1}(V) \cdots \delta_D^{m_n}(V) f^{(n+m)}(tD)\big) + O(t^{N+1-p})
    \]
    for $a \in \A \cup [D,\A]$, $V \in \mB$ bounded and self-adjoint, and $f \in T^{\beta}(\R)$, $\beta< -p$.
\end{rem}

In~\cite[Chapter~5]{EcksteinIochum2018} the question is asked when, for a spectral triple $(\A, \H, D)$, the existence of an asymptotic expansion of $\Tr(e^{-t|D|})$ implies the existence of an expansion of $\Tr(e^{-t|D+V|})$ for a suitable perturbation $V$. Theorem~\ref{T:AsympExpNCG} compared with Example~\ref{E:Counterexample} suggests that this is not generally possible. We illustrate this with the following example.

\begin{ex} 
    Let us revisit Example~\ref{E:Counterexample} where $(\A, \H, D) = (\ell_\infty(\Z_{\geq 1}), \ell_2(\Z_{\geq 1}), D)$, and $D$ is defined by
    \[
        De_n = ne_n,\quad n\geq 1.
    \]
    We take as before
    \[
    a:= \sum_{n=2}^\infty \frac{1}{\log n} e_n.
    \]
    Noting that $|D+a| - |D| = a$, we can apply Theorem~\ref{T:AsympExpNCG} to get (in this situation we can choose $\varepsilon = 0$)
    \[
    \Tr(e^{-t|D+a|}) = \Tr(e^{-tD}) -t \Tr(ae^{-tD}) + O(t).
    \]
    Taking the Mellin transform we get for $\Re(s)>1$ (see~\cite[Proposition~2.10]{EcksteinIochum2018})
    \begin{align*}
    \Tr(|D+a|^{-s}) &= \frac{1}{\Gamma(s)}\int_0^\infty t^{s-1}\Tr(e^{-t|D+a|})\,dt\\
        &= \frac{1}{\Gamma(s)}\int_0^1 t^{s-1}\Tr(e^{-t|D+a|})\,dt + \frac{1}{\Gamma(s)}\int_1^\infty t^{s-1}\Tr(e^{-t|D+a|})\,dt.
\end{align*}
Since
\[
    \Tr(e^{-t|D+a|)}) \leq \Tr(e^{-tD}) = (e^t-1)^{-1} \leq 2 e^{-t},\quad t\geq 1,
\]
we have that
\[
    s\mapsto \frac{1}{\Gamma(s)}\int_1^\infty t^{s-1}\Tr(e^{-t|D+a|})\,ds
\]
is holomorphic.
It follows that
\begin{align*}
    \Tr(|D+a|^{-s}) &= \frac{1}{\Gamma(s)}\int_0^1 t^{s-1}\Tr(e^{-t|D+a|})\,dt + \mathrm{holo}_{\C}(s)\\
    &=   \frac{1}{\Gamma(s)}\int_0^1 t^{s-1} \Tr(e^{-tD}) -t \Tr(ae^{-tD}) dt + \mathrm{holo}_{\Re(s)>-1}(s)\\
    &= \zeta_{0,D}(s) - s \zeta_{a,D}(s+1) + \mathrm{holo}_{\Re(s)>-1}(s).
\end{align*}
Since $s \mapsto \zeta_{a,D}(s+1)$ does not extend holomorphically to any punctured neighbourhood of $s=0$ as the computations in Example~\ref{E:Counterexample} show, we conclude that $\zeta_{a,D}(s) = \Tr(|D+a|^{-s})$ does not admit a meromorphic extension to the entire complex plane. By Theorem~\ref{T:ExpansionToZeta}, 
\[
\Tr (e^{-t|D+a|})
\]
does not admit an asymptotic expansion of the type listed in the theorem.
\end{ex}

We conclude by remarking that the existence of an asymptotic expansion of  $\Tr( Q e^{-t D^2})$ for $Q \in \mB$ of the type in Theorem~\ref{T:ExpansionToZeta} is guaranteed for commutative spectral triples~\cite[Theorem~2.7]{GrubbSeeley1995}, and by the same theorem also for almost commutative spectral triples~\cite[Chapter~8]{Suijlekom2015}.

\newpage

\begin{appendices}
\section{Adjoints}
\label{S:AppEllipticAdj}
In this appendix we discuss various ways of defining the adjoint in $\op^r(\Theta)$. 
\begin{defn}
    Let $A \in \op^r(\Theta)$ so that $A$ extends to a bounded operator
    \[
    A: \H^{s+r} \to \H^{s}
    \]
for all $s\in \mathbb{R}$.
\begin{enumerate}
    \item The adjoint of $A$ as an endomorphism of the topological vector space $\H^\infty$ we denote 
\[
A^\dag : \H^{-\infty} \to \H^{-\infty}
\]
defined by the identity
\[
\langle A u, v \rangle_{(\H^\infty, \H^{-\infty})} = \langle u, A^\dag v \rangle_{(\H^{\infty}, \H^{-\infty})}, \quad u \in \H^\infty, v\in \H^{-\infty}.
\]
\item In similar fashion we denote the adjoint
\[
A'^{_s}: \H^{-s} \to \H^{-s-r}
\]
defined by the relevant identity
\[
\langle Au, v\rangle_{(\H^s, \H^{-s})} = \langle u, A'^{_s}v \rangle_{(\H^{s+r}, \H^{-s-r})}, \quad u \in \H^{s+r}, v\in \H^{-s}.
\]
\item We define the Hermitian adjoint
\[
A^{\flat_s}: \H^{s} \to \H^{s+r}
\]
via the identity
\[
\langle Au, v \rangle_{\H^s} = \langle u, A^{\flat_s} v \rangle_{\H^{s+r}}, \quad u \in \H^{s+r}, v\in \H^s.
\]
\item In case $r \geq 0$, the map
\[
A: \H^{s+r} \subseteq \H^s \to \H^s
\]
  is an unbounded operator on the Hilbert space $\H^s$, so we define another Hermitian adjoint
\[
A^{*_s}: \mathcal{D}_s \to \H^{s},
\]
with domain
\[
\mathcal{D}_s := \{u \in \H^s \ |\  \exists v \in \H^s \forall \phi \in \H^{s+r} : \langle u, T\phi \rangle_{\H^s} = \langle v , \phi \rangle_{\H^{s}} \},
\]
such that
\[
\langle Au, v \rangle_{\H^s} = \langle u, A^{*_s} v \rangle_{\H^{s}}, \quad u \in \H^{s+r}, v\in \mathcal{D}_s.
\]
\end{enumerate}
\end{defn}

These adjoints are related in the following way.
\begin{prop}\label{P:opadjoints}
    Let $A \in \op^r(\Theta)$. Then, for all $s \in \R$,
    \begin{enumerate}
        \item $A'^{_s} = A^\dag \big|_{\H^{-s}}$;
        \item $A^{\flat_s}= \Theta^{-2s-2r} A^\dag \Theta^{2s}\big|_{\H^s}$.
    \end{enumerate}
    If $r\geq 0$, 
    \begin{enumerate}
        \item[3.] $A^{*_s} = \Theta^{-2s} A^\dag \Theta^{2s}\big|_{\mathcal{D}_s}$.
    \end{enumerate}
\end{prop}
\begin{proof}
    \begin{enumerate}
        \item Take $u\in \H^{\infty}\subseteq \H^{s+r}$ and $v\in \H^{-s} \subseteq \H^{-\infty}$. Then
        \begin{align*}
           \langle Au, v \rangle_{(\H^s, \H^{-s})} &= \langle u, A'^{_s} v\rangle_{(\H^{s+r}, \H^{-s-r})}\\
           &= \langle u, A'^{_s} v \rangle_{(\H^{\infty}, \H^{-\infty})}.
        \end{align*}
        We also have
        \begin{align*}
            \langle Au, v \rangle_{(\H^s, \H^{-s})} &  = \langle Au, v \rangle_{(\H^\infty, \H^{-\infty})}\\
            &= \langle u, A^\dag v \rangle_{(\H^{\infty}, \H^{-\infty})}.
        \end{align*}
        Hence it follows that
        \[
        A'^{_s} v = A^\dag v \in \H^{-\infty}, \quad v \in \H^{-s}.
        \]
        
        \item Take $u \in \H^{\infty}$, $v\in \H^{s}$. Then on the one hand,
        \begin{align*}
            \langle Au, v \rangle_{\H^s} &= \langle Au, \Theta^{2s} v \rangle_{(\H^\infty, \H^{-\infty})}\\
            &= \langle u, A^\dag \Theta^{2s} v\rangle_{(\H^{\infty}, \H^{-\infty})},
        \end{align*}
        and on the other hand
        \begin{align*}
            \langle Au, v \rangle_{\H^s} &= \langle u, A^{\flat_s} v \rangle_{\H^{s+r}}\\
            &= \langle u, \Theta^{2s+2r} A^{\flat_s} v \rangle_{(\H^{\infty}, \H^{-\infty})}.
        \end{align*}
        We therefore find
        \[
        A^{\flat_s} = \Theta^{-2s-2r} A^\dag \Theta^{2s}\big|_{\H^s}.
        \]
        \item Take $u \in \H^\infty$, $v \in \mathcal{D}_s \subseteq \H^s \subseteq \H^{-\infty}$. Then
        \begin{align*}
            \langle Au, v\rangle_{\H^s} &= \langle Au, \Theta^{2s} v \rangle_{(\H^\infty, \H^{-\infty})}\\
            &= \langle u, A^\dag \Theta^{2s} v \rangle_{(\H^\infty, \H^{-\infty})},
        \end{align*}
        and
        \begin{align*}
            \langle Au, v\rangle_{\H^s} &= \langle u, A^{*_s} v\rangle_{\H^s}\\
            & =\langle u, \Theta^{2s} A^{*_s} v \rangle_{(\H^\infty, \H^{-\infty})}.
        \end{align*}
        Hence
        \[
        A^{*_s} = \Theta^{-2s} A^\dag \Theta^{2s} \big|_{\mathcal{D}_s}.
        \]
    \end{enumerate}
\end{proof}

An important takeaway from this proposition is that if $A : \H^\infty \to \H^\infty$, we have a priori that
\[
A^\dag: \H^{-\infty} \to \H^{-\infty},
\]
but if $A \in \op^r(\Theta)$ we have in fact that $A^\dag \in \op^r(\Theta)$ (or, more precisely, $A^\dag \big|_{\H^\infty} \in \op^r(\Theta)$). 

It is now also clear that the Hermitian adjoints $A^{\flat_s}$ and $A^{*_s}$ in general cannot be regarded as operators in $\op(\Theta)$ as the operators $A^{\flat_s}$ and $A^{\flat_t}$ do not agree on the intersection $\H^s \cap \H^t$ for $s\not = t$, and the same holds for $A^{*_s}$.

\begin{prop}\label{P:symmetric}
    If $A \in \op^r(\Theta)$, $r \geq 0$, then 
    \[
 A: \H^r\subseteq \H^0 \to \H^0
 \] is symmetric if and only if $A = A^\dag$. 
\end{prop}
\begin{proof}
Suppose that $A = A^\dag$. Let $u\in \H^\infty, v \in \H^{r}$. Then
    \begin{align*}
        \langle Au, v \rangle_{\H^0} &= \langle Au, v\rangle_{(\H^\infty, \H^{-\infty})}\\
        &= \langle u, A^\dag v \rangle_{(\H^\infty, \H^{-\infty})}\\
        &= \langle u, Av \rangle_{\H^0}.
    \end{align*}
 By density of $\H^\infty \subseteq \H^r$, the above equality holds for $u \in \H^r$ as well, and hence $A: \H^r\subseteq \H^0 \to \H^0$ is symmetric.

 On the other hand, if $ A: \H^r\subseteq \H^0 \to \H^0$
 is symmetric, then for $u,v \in \H^\infty,$ 
 \begin{align*}
     \langle u, A^\dag v \rangle_{(\H^\infty, \H^{-\infty})} &= \langle Au, v\rangle_{\H^0} \\
    &= \langle u, Av \rangle_{\H^0}\\
    &= \langle u, Av \rangle_{(\H^\infty, \H^{-\infty})},
 \end{align*}
 showing that $A^\dag v = Av \in \H^\infty$ which implies that $A=A^\dag \in \op^r(\Theta)$.
\end{proof}

\section{Functional Calculus for \texorpdfstring{$\op^0(\Theta)$}{op0}}\label{S:AppFunctCalc0}
In Section~\ref{S:FunctCalc} we proved that symmetric $\Theta$-elliptic operators in $\op^r(\Theta)$ for $r>0$ admit a functional calculus. The approach of that section does not apply for the case $r=0$. To illustrate how different the zero-order case is, consider the situation where $\Theta = (1+\Delta)^{1/2}$ on $L_2(M)$ where $M$ is a compact subset of $\R^d$. We have that for $\phi: M \to \mathbb{R}$, the multiplication operator
    \begin{align*}
        M_\phi: L_2(M) &\to L_2(M)\\
        \xi &\mapsto \phi \cdot \xi,
    \end{align*}
    where $\phi \cdot \xi (x) = \phi(x) \xi(x)$, can only be in $\op(\Theta)$ if $\phi$ is smooth. For smooth $\phi$ we have that $M_\phi \in \op^0(\Theta)$. If $f(M_\phi) \in \op(\Theta)$ for all $M_\phi \in \op^0(\Theta)$, then the identity 
    \[
    f(M_\phi) = M_{f \circ \phi}
    \]
    shows that the function $f$ has to be smooth itself and no functional calculus with general functions in $L^\beta_\infty(E)$ is possible.

For $\op^0(\Theta)$ we therefore use a different strategy altogether. 
An approach by Davies~\cite{Davies1995c, Davies1995a} on the construction of a functional calculus using almost analytic extensions directly applies. Using almost analytic extensions to obtain a functional calculus for pseudodifferential operators has precedent in the works of, amongst others,  H\"ormander~\cite{Hormander1969}, Helffer--Sj\"ostrand~\cite{HelfferSjostrand1989}, Dimassi--Sj\"ostrand~\cite[Chapter 8]{DimassiSjostrand1999}, and Bony~\cite{Bony2013}.

Recall the definition of an almost analytic extension, Definition~\ref{def:AlmostAnalyticExtension1}.

\begin{thm}[\cite{Davies1995a}]\label{T:Helffer-Sjostrand}
    Let $f \in C^\infty_c(\mathbb{R})$ with almost analytic extension $\tilde{f}$ as in Definition~\ref{def:AlmostAnalyticExtension1}, so that
    \[
    f(x) = -\frac{1}{\pi} \int_\mathbb{C} \left( \frac{\partial \tilde{f}}{\partial \overline{z}}(z)\right) (z-x)^{-1} dz, \quad x \in \mathbb{R}.
    \]
    For any closed, densely defined operator $H$ with $\sigma(H) \subseteq \mathbb{R}$, if for some $\alpha \in \mathbb{R}_{\geq 0}$ we have the estimate
    \[
    \| (z-H)^{-1} \| \leq C \frac{1}{|\Im(z)|}\left(\frac{\langle z \rangle}{|\Im(z)|}\right)^\alpha, \quad z\in \mathbb{C}\subseteq \mathbb{R},
    \]
    then we have that
    \[
    f(H) := -\frac{1}{\pi} \int_\mathbb{C} \frac{\partial \tilde{f}}{\partial \overline{z}} (z-H)^{-1} dz
    \]
    defines a bounded operator on $\H$ independent of the choice of $N > \alpha$ and $\tau$ in the construction of the extension $\tilde{f}$, with
    \[
    \|f(H)\|_\infty \leq \sum_{k=0}^{N+1}\|f\|_{T^0(\R), k}.
    \]
    The integral should be interpreted as a $B(\H)$-valued Bochner integral. In case $H$ is self-adjoint, this agrees with the continuous functional calculus.
\end{thm}
We thank Dmitriy Zanin for providing a key step in the following proof, which is an adaptation of an argument by Beals~\cite[Lemma~3.1]{Beals1977}.

\begin{prop}\label{P:Invertop}
    Let $X \in \op^r(\Theta)$ be such that $[\Theta, X] \in \op^r(\Theta)$. If the extension 
    \[
    X: \H^{s_0+r} \to \H^{s_0}
    \]
    has a bounded inverse 
    \[
    X^{-1}:\H^{s_0} \to \H^{s_0+r}
    \]
    for one particular $s_0 \in \mathbb{R}$, then $X^{-1}\big|_{\H^\infty} \in \op^{-r}(\Theta)$. We have $X X^{-1}|_{\H^\infty} = X^{-1}X|_{\H^{\infty}} = 1_{\H^\infty}$. In particular, if $X \in \op^r(\Theta)$ and $[\Theta, X] \in \op^r(\Theta)$ with $r \geq 0$, then we have as (unbounded) operators 
    \[
    \sigma(X: \H^{s_0 +r} \subseteq \H^{s_0} \to \H^{s_0}) = \sigma(X: \H^{s+r}\subseteq \H^s \to \H^s)
    \]
    for all $s\in \mathbb{R}$, where $\sigma$ denotes the spectrum of the operator.
\end{prop}
\begin{proof}
Since $\op^r(\Theta) = \op^0(\Theta) \cdot \Theta^r$, it suffices to prove the proposition for $r=0$. 

Suppose that $X \in \op^0(\Theta)$ is a bijection on $\H^{s_0} \to \H^{s_0}$ and write $X^{-1}:\H^{s_0} \to \H^{s_0}$. Then $X$ restricts to a necessarily injective map $\H^{{s_0}+1} \to \H^{{s_0}+1}$. We now prove that $X: \H^{{s_0}+1} \to \H^{{s_0}+1}$ is also surjective.

We follow~\cite[Lemma~3.1]{Beals1977}, filling in some omitted details. Take $v \in \H^{s_0 + 1}$, then there exists $u \in \H^{s_0}$ with $Xu = v$. Let $\varepsilon >0$, then $\frac{\Theta}{1+\varepsilon \Theta} \in \op^0(\Theta)$, and
\begin{align*}
    \frac{\Theta}{1+ \varepsilon \Theta} u &= \frac{\Theta}{1+ \varepsilon \Theta} X^{-1} v \\
    &= X^{-1} \frac{\Theta}{1+\varepsilon \Theta}v + X^{-1}\left[\frac{\Theta}{1+ \varepsilon \Theta}, X\right] X^{-1}v.
\end{align*}
Now, the $\H^{s_0}$ norm of the right-hand side is bounded independent of $\varepsilon$: 
\begin{align*}
    \left \| \frac{\Theta}{1+\varepsilon \Theta} \right \|_{\H^{s_0+1} \to \H^{s_0}} &\leq \| \Theta \|_{\H^{s_0+1} \to \H^{s_0}};\\
\left\|\left[\frac{\Theta}{1+ \varepsilon \Theta}, X\right] \right\|_{\H^{s_0} \to \H^{s_0}} &= \|(1+\varepsilon\Theta)^{-1}[\Theta,X](1+\varepsilon \Theta)^{-1}\|_{\H^{s_0}\to H^{s_0}} \leq \|[\Theta, X]\|_{\H^{s_0}\to \H^{s_0}}.
\end{align*}
This implies that $u \in \H^{s_0+1}$, a fact that can be quickly verified with the spectral theorem and Fatou's lemma. Therefore, 
\[
    X: \H^{s_0 + 1} \to \H^{s_0+1}
\]
is a bijection. By induction and interpolation (Proposition~\ref{P:Interpolation}), the same assertion holds for each $\H^s$, $s \geq s_0$.

Finally, it is a basic fact that the adjoint of a bijective operator is bijective, i.e.
\[
X'^{_{s_0}} = X^\dag|_{\H^{-s_0}}:\H^{-s_0} \to \H^{-s_0}
\]
is a bijection. Since
\[
[\Theta, X^\dag] = - [\Theta, X]^{\dag} \in \op^0(\Theta),
\]
we can apply the same arguments as above to deduce that
\[
X^{\dag}:\H^{-s} \to \H^{-s}
\]
is a bijection for all $-s \geq -s_0$. This implies that 
\[
X = X^{\dag \dag}: \H^{s} \to \H^{s}
\]
is a bijection for all $s \leq s_0$.
\end{proof}
This following type of estimate on the resolvent also appears in $L_p$-boundedness problems, see \cite{Davies1995b, Davies1995a, JensenNakamura1994}.
\begin{lem}\label{L:ResolventEstimate}
    Let $A \in \op^0(\Theta)$ be such that $[\Theta, A] \in \op^0(\Theta)$ and $\overline{A}^{0,0}:\H \to \H$ is self-adjoint. 
    Then for all $s \in \mathbb{R},$ there is a constant $C_s > 0$ such that
    \begin{align*}
        \|(z-A)^{-1}\|_{\H^s\to \H^s} &\leq C_s \frac{1}{|\Im(z)|} \left(\frac{\langle z \rangle}{|\Im(z)|}\right)^{2^{|s|}-1}, \quad z \in \mathbb{C}\setminus \mathbb{R}.
    \end{align*}
\end{lem}
\begin{proof}
    The proof is by induction and interpolation (Proposition~\ref{P:Interpolation}). For $s=0$, the estimate holds by self-adjointness of $\overline{A}^{0,0}$. Note that $(z-A)^{-1} \in \op^0(\Theta)$ due to Proposition~\ref{P:Invertop}.
    
    Suppose the inequality is proved for a fixed $s\in\R_{\geq 0}$. Then for $z \in \mathbb{C}\setminus \mathbb{R}$,
    \begin{align*}
        & \|(z-A)^{-1}\|_{\H^{s+1}\to \H^{s+1}} \\
        &= \|\Theta(z-A)^{-1}\Theta^{-1}\|_{\H^s\to \H^s}\\
        &\leq \|(z-A)^{-1}\|_{\H^s\to \H^s}+\|(z-A)^{-1}[\Theta,A](z-A)^{-1}\Theta^{-1}\|_{\H^s\to \H^s}\\
        &\leq  \|(z-A)^{-1}\|_{\H^s\to \H^s} \Big(1+ \|[\Theta,A](i+A)^{-1}\|_{\H^s\to \H^s}\|(i+A)(z-A)^{-1}\|_{\H^s\to \H^s} \|\Theta^{-1}\|_{\H^s \to \H^s} \Big).
    \end{align*}
    Note that $(i+A)^{-1} \in \op^{0}(\Theta)$ by Proposition~\ref{P:Invertop}, so that for some constant $B_s >0$,
    \[
    \|[\Theta,A](i+A)^{-1}\|_{\H^s\to \H^s} \|\Theta^{-1}\|_{\H^s\to \H^s}\leq B_s.
    \]
    Using the resolvent identity, we have
    \[
        (i+A)(z-A)^{-1} = (i+z)(z-A)^{-1}-1
    \]
    and therefore
    \begin{align*}
        \|(i+A)(z-A)^{-1}\|_{\H^s\to \H^s} & \leq 1+ |z+i|\|(z-A)^{-1}\|_{\H^s\to \H^s}.
    \end{align*}
    This yields
    \begin{align*}
        \|(z-A)^{-1}\|_{\H^{s+1}\to \H^{s+1}} &\leq  \|(z-A)^{-1}\|_{\H^s\to \H^s} (1+B_s )\\
        &\quad + |z+i| B_s  \|(z-A)^{-1}\|_{\H^s\to \H^s}^2 \\
        &\leq (1+B_s) \|(z-A)^{-1}\|_{\H^s\to \H^s} \cdot (1+ |z+i|\|(z-A)^{-1}\|_{\H^s\to \H^s}).
    \end{align*}
    This estimate also holds with $|z-i|$ on the right-hand side, and $\min(|z+i|,|z-i|)\leq\langle z\rangle$, so that
    \begin{align*}
        \|(z-A)^{-1}\|_{\H^{s+1}\to \H^{s+1}} &\leq  (1+B_s)  \|(z-A)^{-1}\|_{\H^s\to \H^s} \cdot (1+ \langle z \rangle \|(z-A)^{-1}\|_{\H^s\to \H^s}),
    \end{align*}
    from which the claimed estimate follows. Induction and interpolation now provide the estimate for all $s \geq 0$. 

    The case $s \leq 0$ is proved in the same manner, using induction in the negative direction. Namely, the norm 
    \[
    \| (z-A)^{-1}\|_{\H^{s-1} \to \H^{s-1}} = \|\Theta^{-1}(z-A)^{-1}\Theta \|_{\H^s\to \H^s}
    \]
    can be estimated as before.
\end{proof}

\begin{prop}\label{P:FunctCalcOp0}
    Let $A \in \op^0(\Theta)$ be such that $[\Theta, A] \in \op^0(\Theta)$ and $\overline{A}^{0,0}: \H \to \H$ is self-adjoint. If $f \in C^\infty(\R)$, then $f(A) \in \op^{0}(\Theta)$ with
    \[
    \|f(A)\|_{\H^s \to \H^s} \leq \sum_{k=0}^{\lceil 2^{|s|}\rceil+1} \|f\|_{T^0(\R), k}.
    \]
\end{prop}
\begin{proof}
Without loss of generality, we assume that $f\in C_c^\infty(\R)$. As $A\in \op^0(\Theta)$, it extends to a bounded operator
\[
\overline{A}^{s,s}: \H^s \to \H^s, \quad s\in \R.
\]
Furthermore, by Proposition~\ref{P:InvertElliptic}, we have
\[
(A-z)^{-1} \in \op^0(\Theta), \quad z \in \C\setminus\R.
\]
Theorem~\ref{T:Helffer-Sjostrand} and Lemma~\ref{L:ResolventEstimate} combined give that
\[
f(\overline{A}^{s,s}): \H^s\to \H^s, \quad s \in \R,
\]
is a bounded operator with the norm bound as claimed. By construction, for $\xi \in \H^\infty$ we have that
\begin{align*}
    f(\overline{A}^{s,s})\xi = \int_{\C} \frac{\partial \tilde{f}}{\partial \overline{z}} (\overline{A}^{s,s} - z)^{-1} \xi dz \in \H^s,
\end{align*}
as an $\H^s$-valued Bochner integral. It is clear that
\[
(\overline{A}^{s,s} - z)^{-1} \big|_{\H^\infty} = (A-z)^{-1},
\]
and therefore these arguments show that for $\xi \in \H^\infty$ the integral
\[
\int_{\C} \frac{\partial \tilde{f}}{\partial \overline{z}} (A - z)^{-1} \xi dz
\]
can be evaluated as a Bochner integral in each Sobolev space $\H^s$. Hence
\[
\xi \mapsto \int_{\C} \frac{\partial \tilde{f}}{\partial \overline{z}} (A - z)^{-1} \xi dz
\]
forms a bounded linear map on $\H^\infty$; denote this operator by $f(A): \H^\infty \to \H^\infty$. Then since $f(A)$ agrees with $f(\overline{A}^{s,s})$ on $\H^\infty$, we must have
\[
\overline{f(A)}^{s,s} = f(\overline{A}^{s,s}),
\]
and thus we have $f(A) \in \op^0(\Theta)$.
\end{proof}

This functional calculus can be used to construct multiple operator integrals, similarly to Theorem~\ref{T:MainMOIConstruction}.
\begin{thm}
    Let $H_0, \ldots, H_n, [\Theta, H_0], \ldots, [\Theta, H_n] \in \op^0(\Theta)$ be such that each $\overline{H_i}^{0,0}$ is self-adjoint. For $\phi \in T^0(\R) \boxtimes_i \cdots \boxtimes_i T^0(\R)$ (see Definition~\ref{def:ProjIntBox}) with corresponding representation
    \[
\phi(\lambda_0,\ldots,\lambda_n)=\int_\Omega a_0(\lambda_0,\omega)\cdots a_n(\lambda_n,\omega) d\nu(\omega),
\]
the integral
\[
T_\phi^{H_0,\ldots,H_n}(X_1,\ldots,X_n)\psi:=\int_\Omega a_0(H_0,\omega)X_1 a_1(H_1,\omega)\cdots X_n a_n(H_n,\omega)\psi\,d\nu(\omega), \quad \psi\in\H^\infty
\]
for $X_1,\ldots,X_n\in\op(\Theta)$, converges as a Bochner integral in $\H^s$ for every $s \in \R$, and defines an $n$-multilinear map $T_\phi^{H_0,\ldots,H_n}:\op^{r_1}(\Theta)\times\cdots\times\op^{r_n}(\Theta)\to\op^{\sum_j r_j}(\Theta)$.
\end{thm}
\begin{proof}
    Recall that by definition of $T^0(\R) \boxtimes_i \cdots \boxtimes_i T^0(\R)$ we have for all $k_0, \ldots, k_n \geq 0$
    \[
    \int_\Omega \|a_0(\cdot, \omega)\|_{T^0(\R), k_0} \cdots \|a_n(\cdot, \omega)\|_{T^0(\R), k_n} d|\nu|(\omega) < \infty.
    \]
    The proof of the theorem is then identical to the proof of Theorem~\ref{T:MainMOIConstruction}, using that
    \[
    \|a_j(H_j, \omega)\|_{\H^s \to \H^s} \leq \sum_{k=0}^{\lceil 2^{|s|}\rceil+1} \|a_j(\cdot, \omega)\|_{T^0(\R), k}
    \]
    instead of
    \[
        \|a_j(H_j,\omega)\|_{\H \indices{^{s+k_j}}\to\H\indices{^s}}\leq C_{s,H_j}\|a_j(\cdot,\omega)\|_{L^{\beta_j}_{\infty}(E_j)}.
    \]
\end{proof}
Due to Lemma~\ref{L:DivDifT}, the divided difference $f^{[n]}$ for $f \in T^\alpha(\R)$ with $\alpha < n$ is in particular a permitted symbol.

One may wonder how the condition $A, [\Theta,A] \in \op^r(\Theta)$ compares to the condition of $A$ being $\Theta$-elliptic with $r >0$. 
In principle, these conditions are incomparable:
\begin{itemize}
    \item On $\R^d$, the operator $\partial_{x_1} \in \op^1((1-\Delta)^{\frac12})$ satisfies the first condition but is not $(1-\Delta)^{\frac12}$-elliptic.
    \item Taking $\Theta=1-\partial_x^2$ on $\R$, the dilation $T_\varepsilon f(x):=f(\varepsilon x)$ by $\varepsilon\neq 0$ is in $\op^0(1-\partial_x^2)$ with inverse $T_{\varepsilon^{-1}}$. Hence, $(1-\partial_x^2) \circ T_\varepsilon \in \op^1(1-\partial_x^2)$ is $(1-\partial_x^2)$-elliptic (and even invertible on the nose). However, as once can check by an elementary computation, 
    \[
    [(1-\partial_x^2), (1-\partial_x^2) \circ T_\varepsilon] = T_\varepsilon \circ \big((\varepsilon^4-\varepsilon^2)\partial_x^4 + (1-\varepsilon^2)\partial_x^2\big) \not\in \op^1(1-\partial_x^2),
    \] 
    if $\varepsilon \neq \pm 1$.
\end{itemize}

The following consequence of Proposition~\ref{P:Invertop} shows that, if $A$ has a non-empty resolvent set on the domain $\H^r$, then $\Theta$-ellipticity follows from the condition that $A, [\Theta,A] \in \op^r(\Theta)$. The second example above shows that under this additional assumption, $\Theta$-ellipticity is therefore strictly weaker condition: $(1+i \partial_x)\circ T_\varepsilon$ has non-empty resolvent set since it is invertible, but does not satisfy the commutator condition.
\begin{cor}\label{C:EllipticEquiv}
    Let $A\in \op^r(\Theta), r>0,$ be such that
    \[
    A:\H^r\subseteq \H^0 \to \H^0
    \]
    has a non-empty resolvent set (for example if $A$ is self-adjoint with domain $\H^r$), and suppose that $[\Theta, A] \in \op^r(\Theta)$ as well. Then $A$ is $\Theta$-elliptic.
\end{cor}
\begin{proof}
    By assumption, there exists $z \in \C$ such that 
    \[
    z-A: \H^r \subseteq \H^0 \to \H^0
    \]
    is invertible, and because $[\Theta, A] \in \op^r(\Theta)$ it follows that $(z-A)^{-1} \in \op^{-r}(\Theta)$ by Proposition~\ref{P:Invertop}. Now,
    \[
    A (z-A)^{-1} = - 1 + z (z-A)^{-1}. 
    \]
    In other words, $-(z-A)^{-1}$ is an inverse of $A$ modulo $\op^{-r}(\Theta)$. Corollary~\ref{C:EllipticParametrix} gives that $A$ is $\Theta$-elliptic.
\end{proof}

\section{Combinatorial Computations}
\label{S:AppCombinComp}
This appendix is dedicated to proving Proposition~\ref{P:Expansion}. The computations are standard, cf.~\cite[Lemma~6.11]{CPRS1}, \cite[Equation~(71)]{ConnesMoscovici1995}, \cite[Lemma~2.12]{Higson2003} -- the novelty is that they can be performed in the very general context of (unbounded) MOIs.

\begin{defn}
    The multiset coefficient $\multinom{n}{k}$ for $n, k \in \mathbb{Z}_{\geq 0}$ is defined as 
    \[
    \multinom{n}{k} := \binom{n+k-1}{k}.
    \]
\end{defn}

\begin{lem}\label{L:Commute1} For $f\in T^\beta(\R)$, $H\in \op^h(\Theta)$, $h> 0$ symmetric and $\Theta$-elliptic, $X_i \in \op^{r_i}(\Theta)$, we have
\begin{align*}
    &T^{H} _{f^{[n+j]}}(\underbrace{1, \dots, 1}_{j}, X_1, \dots, X_n ) \\
    &= \sum_{m=0}^{N}  \multinom{m+1}{j} \delta_H^{m}(X_1) T^{H}_{f^{[n+j+m]}}(\underbrace{1, \dots, 1}_{j+1+m}, X_2, \dots, X_n) + R^n_{j,N}(X_1, \dots, X_n),
\end{align*}
where 
\[
R^n_{j,N}(X_1, \dots, X_n) :=  \sum_{l=0}^{j} \multinom{N+1}{l} T^{H}_{f^{[n+j+N+1]}}(\underbrace{1, \dots, 1}_{j-l}, \delta_H^{N+1}(X_1), \underbrace{1, \dots, 1}_{N+1+l}, X_2, \dots, X_n).
\]
\end{lem}
\begin{proof}
    Multiset coefficients have the property that
    \[
    \sum_{l=0}^j \multinom{m}{l} = \multinom{m+1}{j}.
    \]
    The assertion of the lemma follows by induction on $N$. For $N=0$,
    \begin{align*}
        &T^{H}_{f^{[n+j]}}(\underbrace{1, \dots, 1}_{j}, X_1, \dots, X_n )\\
        &= X_1 T^{H}_{f^{[n+j]}}(\underbrace{1, \dots, 1}_{j+1}, X_2, \dots, X_n ) + \sum_{l=0}^j T^{H}_{f^{[n+j+1]}}(\underbrace{1, \dots, 1}_{j-l}, \delta_H(X_1), \underbrace{1, \dots, 1}_{1+l}, X_2, \dots, X_n)
    \end{align*}
    by applying Proposition~\ref{P:UMOIcom} $j+1$ times on $X_1$.

    Suppose that the assertion holds for $N-1$. Then 
    \begin{align*}
         &T^{H}_{f^{[n+j]}}(\underbrace{1, \dots, 1}_{j}, X_1, \dots, X_n ) \\
         &= \sum_{m=0}^{N-1}  \multinom{m+1}{j} \delta_H^m(X_1) T^{H}_{f^{[n+j+m]}}(\underbrace{1, \dots, 1}_{j+1+m}, X_2, \dots, X_n) \\
    & \quad +  \sum_{l=0}^{j} \multinom{N}{l} T^{H}_{f^{[n+j+N]}}(\underbrace{1, \dots, 1}_{j-l}, \delta_H^{N}(X_1), \underbrace{1, \dots, 1}_{N+l}, X_2, \dots, X_n)\\
    &=^*  \sum_{m=0}^{N-1} \multinom{m+1}{j} \delta_H^m(X_1) T^{H}_{f^{[n+j+m]}}(\underbrace{1, \dots, 1}_{j+1+m}, X_2, \dots, X_n) \\
    & \quad +  \sum_{l=0}^{j} \multinom{N}{l} \delta_H^{N}(X_1) T^{H}_{f^{[n+j+N]}}(\underbrace{1, \dots, 1}_{N+1+j}, X_2, \dots, X_n)\\
    & \quad +  \sum_{l=0}^{j} \multinom{N}{l} \sum_{k=0}^{j-l} T^{H}_{f^{[n+j+N+1]}}(\underbrace{1, \dots, 1}_{j-l-k}, \delta_H^{N+1}(X_1), \underbrace{1, \dots, 1}_{N+1+l+k}, X_2, \dots, X_n),
    \end{align*}
    where in the step marked with $*$ we applied Proposition~\ref{P:UMOIcom} $j-l$ times on $\delta_H^{N}(X_1)$. Continuing on, 
    \begin{align*}
 &T^{H}_{f^{[n+j]}}(\underbrace{1, \dots, 1}_{j}, X_1, \dots, X_n )\\
    &=  \sum_{m=0}^{N-1} \multinom{m+1}{j} \delta_H^m(X_1) T^{H}_{f^{[n+j+m]}}(\underbrace{1, \dots, 1}_{j+1+m}, X_2, \dots, X_n) \\
    & \quad +   \multinom{N+1}{j} \delta_H^{N}(X_1) T^{H}_{f^{[n+j+N]}}(\underbrace{1, \dots, 1}_{N+1+j}, X_2, \dots, X_n)\\
    & \quad +  \sum_{l=0}^{j} \sum_{k=0}^{j-l}  \multinom{N}{l} T^{H}_{f^{[n+j+N+1]}}(\underbrace{1, \dots, 1}_{j-l-k}, \delta_H^{N+1}(X_1), \underbrace{1, \dots, 1}_{N+1+l+k}, X_2, \dots, X_n).
    \end{align*}
 In the last sum, relabel $r:= k+l$, so that
\begin{align*}
    &T^{H}_{f^{[n+j]}}(\underbrace{1, \dots, 1}_{j}, X_1, \dots, X_n ) \\
    &= \sum_{m=0}^{N}  \multinom{m+1}{j} \delta_H^m(X_1) T^{H}_{f^{[n+j+m]}}(\underbrace{1, \dots, 1}_{j+1+m}, X_2, \dots, X_n) \\
    & \quad +  \sum_{r=0}^{j} \sum_{l=0}^{r}  \multinom{N}{l} T^{H}_{f^{[n+j+N+1]}}(\underbrace{1, \dots, 1}_{j-r}, \delta_H^{N+1}(X_1), \underbrace{1, \dots, 1}_{N+1+r}, X_2, \dots, X_n)\\
    &= \sum_{m=0}^{N} \multinom{m+1}{j} \delta_H^m(X_1) T^{H}_{f^{[n+j+m]}}(\underbrace{1, \dots, 1}_{j+1+m}, X_2, \dots, X_n) \\
    & \quad +  \sum_{r=0}^{j}   \multinom{N+1}{r} T^{H}_{f^{[n+j+N+1]}}(\underbrace{1, \dots, 1}_{j-r}, \delta_H^{N+1}(X_1), \underbrace{1, \dots, 1}_{N+1+r}, X_2, \dots, X_n).
\end{align*}
This concludes the induction step.
\end{proof}

\begin{proof}[Proof of Proposition~\ref{P:Expansion}]
Apply Lemma~\ref{L:Commute1} first to the first entry of $T^{H}_{f^{[n]}}(X_1, \dots, X_n)$,
\begin{align*}
    T^{H}_{f^{[n]}}(X_1, \dots, X_n ) &= \sum_{m_1=0}^{N}  \multinom{m_1+1}{0} \delta_H^{m_1}(X_1) T^{H}_{f^{[n+m_1]}}(\underbrace{1, \dots, 1}_{m_1+1}, X_2, \dots, X_n) + R^{n}_{0,N}(X_1, \dots, X_n).
\end{align*}
Apply Lemma~\ref{L:Commute1} once more, expanding up to order $N-m_1$ instead of $N$.
\begin{align*}
    &T^{H}_{f^{[n]}}(X_1, \dots, X_n ) \\
    &= \sum_{m_1=0}^{N} \sum_{m_2=0}^{N-m_1} \multinom{m_1+1}{0} \multinom{m_2+1}{m_1+1} \delta_H^{m_1}(X_1) \delta_H^{m_2}(X_2) T^{H}_{f^{[n+m_1+m_2]}}(\underbrace{1, \dots, 1}_{2+m_1+m_2}, X_3, \dots, X_n)\\
    & \quad + \sum_{m_1=0}^{N}\multinom{m_1+1}{0}  \delta_H^{m_1}(X_1) R^{n-1}_{1+m_1,N-m_1}(X_2, \dots, X_n) + R^{n}_{0,N}(X_1, \dots, X_n).
\end{align*}
Repeating gives the formula
\begin{align*}
    &T^{H}_{f^{[n]}}(X_1, \dots, X_n )\\
    &= \sum_{m=0}^{N}  \sum_{m_1 + \dots + m_n = m} \prod_{j=1}^n  \multinom{m_j+1}{j-1+m_1+\dots + m_{j-1}}\delta_H^{m_1}(X_1) \cdots \delta_H^{m_n}(X_n) T^{H}_{f^{[n+m]}}(\underbrace{1, \dots, 1}_{n+m})\\
    & \quad + S^n_{N}(X_1, \dots X_n),
\end{align*}
where
\begin{align*}
     S^n_{N}(X_1, \dots X_n)&:=   \sum_{k=0}^{n-1} \sum_{m_1 + \dots + m_k \leq N} \prod_{j=1}^{k}  \multinom{m_j+1}{j-1+m_2+\dots + m_{j-1}} \\
     &\quad \times \delta_H^{m_1}(X_1)\cdots \delta_H^{m_k}(X_{k})  R^{n-k}_{k+m_1+\dots+m_{k},N-m_1-\dots-m_{k}}(X_{k+1}, \dots, X_n).
\end{align*}

The observation that 
\[
\multinom{n}{k} = \multinom{k+1}{n-1},
\]
and the definition
\[
\multinom{n}{k} = \binom{n+k-1}{k},
\]
finishes the proof of the proposition.
\end{proof}

\end{appendices}
\printbibliography

\end{document}